\documentclass[review]{elsarticle}

\usepackage{lineno,hyperref}

\usepackage{amssymb}%
\usepackage{amsmath}%
\usepackage{amsthm}%
\usepackage{flexisym}%

\newtheorem{theorem}{Theorem}[section]%
\newtheorem{lemma}[theorem]{Lemma}%
\newtheorem{corollary}[theorem]{Corollary}

\newtheorem{example}[theorem]{Example}
\newtheorem{remark}[theorem]{Remark}

\newtheorem{definition}[theorem]{Definition}

\DeclareMathOperator\arctanh{arctanh}

\DeclareMathOperator\Li{Li}

\allowdisplaybreaks

\modulolinenumbers[5]

\journal{Journal of Mathematical Analysis and Applications }

\begin{document}

\begin{frontmatter}

\title{On the interplay among hypergeometric functions, complete 
 elliptic integrals, and Fourier--Legendre expansions}

\author[Yorku]{John M. Campbell\corref{mycorrespondingauthor}}
\cortext[mycorrespondingauthor]{Corresponding author}
\ead{jmaxwellcampbell@gmail.com}

\author[Pisau]{Jacopo D'Aurizio}

\author[209West97thStreet]{Jonathan Sondow}

\address[Yorku]{Department of Mathematics and Statistics, York University, Toronto, Ontario, M3J 1P3}

\address[Pisau]{Department of Mathematics, University of Pisa, Italy, Largo Bruno Pontecorvo 5 - 56127}

\address[209West97thStreet]{209 West 97th Street, New York, NY 10025}

\begin{abstract}
 Motivated by our previous work on hypergeometric functions and the parbelos constant, we perform a deeper investigation on the 
 interplay among generalized complete elliptic integrals, Fourier--Legendre (FL) series expansions, and ${}_p F_q$ series.
 We produce new 
 hypergeometric transformations and closed-form evaluations for new series involving harmonic 
 numbers, through the use of the integration method outlined as follows: Letting $\text{{\bf K}}$ denote the complete elliptic 
 integral of the first kind, for a suitable function $g$ we evaluate integrals such as $$ \int_{0}^{1} \text{{\bf K}}\left( \sqrt{x} 
 \right) g(x) \, dx $$ in two different ways: (1) by expanding $\text{{\bf K}}$ as a Maclaurin series, perhaps after a 
 transformation or a change of variable, and then integrating term-by-term; and (2) by expanding $g$ as a shifted FL 
 series, and then integrating term-by-term. Equating the expressions produced by these two approaches often gives us new 
 closed-form evaluations, as in the formulas involving Catalan's constant $G$
\begin{align*}
 \sum _{n = 0}^{\infty } \binom{2 n}{n}^2 \frac{H_{n + \frac{1}{4}} - 
 H_{n-\frac{1}{4}}}{16^{n} } 
 & = \frac{\Gamma^4 \left(\frac{1}{4}\right)}{8 \pi^2}-\frac{4 G}{\pi}, \\ 
 \sum _{m, n \geq 0} \frac{\binom{2 m}{m}^2 \binom{2 n}{n}^2 }{ 16^{m + n} (m+n+1) (2 m+3) } & = 
 \frac{7 \zeta (3) - 4 G}{\pi^2}.
\end{align*}
\end{abstract}


\end{frontmatter}

\linenumbers

\section{Introduction}\label{Introduction}
 In the history of mathematical analysis, there are many strategies for computing infinite series in symbolic form and it remains a very 
 active area of research. In our present article, we introduce a variety of new results on the closed-form evaluation of 
 hypergeometric series and harmonic summations through the use of new techniques that are mainly based on the use of complete 
 elliptic integrals and the theory of Fourier--Legendre (FL) expansions. 

 Inspired in part by our previous work on the evaluation of a ${}_{3}F_{2}(1)$ series related to the parbelos constant 
 \cite{CampbellDAurizioSondow}, which, in turn, came about through the discovery \cite{Campbell2018Ramanujan} of an 
 integration technique for evaluating series involving squared central binomial coefficients and harmonic numbers in terms of 
 $\frac{1}{\pi}$, in this article we apply a related integration method to determine new identities for hypergeometric 
 expressions, as well as new evaluations for binomial-harmonic series.

 We recall that a \emph{hypergeometric series} is an infinite series $\sum_{i=0}^{\infty} c_{i}$ such that there exist 
 polynomials $P$ and $Q$ satisfying 
\begin{equation}\label{20190119724PM1A}
 \frac{c_{i+1}}{c_i} = \frac{P(i)}{Q(i)} 
\end{equation}
 for each $i \in \mathbb{N}_{0}$. If $P$ and $Q$ can be written as $P(i) = (i + a_{1})(i + a_{2})\cdots(i+ 
 a_{p}) $ and $ Q(i) = (i + b_1)(i + b_2) \cdots (i + b_q)(i+1)$, a \emph{generalized 
 hypergeometric function} is an infinite series of the form $\sum_{i=0}^{\infty} c_{i} x^{i}$, and is denoted by
\begin{equation*}
 \sum_{i=0}^{\infty} c_{i} x^{i}
 = {}_{p}F_{q}\!\!\left[ \begin{matrix} 
 a_{1}, a_{2}, \ldots, a_{p} \vspace{1mm}\\ 
 b_{1}, b_{2}, \ldots, b_{q} \end{matrix} \ \Bigg| \ x 
 \right]. 
\end{equation*}
 The complete elliptic integrals of the first and second kinds, which are respectively denoted 
 with $\text{{\bf K}}$ and $\text{{\bf E}}$, 
 may be defined as 
\begin{align*}
 & \text{{\bf K}}(k) = 
 \frac{\pi}{2} \cdot {}_{2}F_{1}\!\!\left[ 
 \begin{matrix}
 \frac{1}{2}, \frac{1}{2} \vspace{1mm}\\
 1
 \end{matrix} \ \Bigg| \ k^2 \right] \\
 & \text{{\bf E}}(k) = 
 \frac{\pi}{2} \cdot {}_{2}F_{1}\!\!\left[ 
 \begin{matrix}
 \frac{1}{2}, -\frac{1}{2} \vspace{1mm}\\
 1
 \end{matrix} \ \Bigg| \ k^2 \right]. 
\end{align*} 
 We are interested in evaluating integrals such as 
\begin{equation}\label{20190119234PM1A}
 \int_{0}^{1} \text{{\bf K}}\left( \sqrt{x} \right) g(x) \, dx 
\end{equation}
 for a suitable function $g$, 
 by expanding $\text{{\bf K}}$ as a Maclaurin series, perhaps 
 after a manipulation of the expression $\text{{\bf K}}$, and integrating term-by-term. However, by replacing $g(x)$ with 
 its shifted Fourier--Legendre series expansion, integrating term-by-term and equating the two resulting series, we 
 often obtain new closed-form evaluations. We illustrate this idea with the example described in Section \ref{20190119229PM2A} 
 below that is
 taken from \cite{CampbellDAurizioSondow}, after a preliminary 
 discussion concerning the basics of FL theory.

 Legendre functions of order $n$ are solutions to Legendre's differential equation $$ (1-x^2) \frac{d^{2} y}{dx^{2}} - 2 x 
 \frac{dy}{dx} + n(n+1) y = 0 $$ for $n > 0$ and $|x| < 1$. For $n \in \mathbb{N}_{0}$, \emph{Legendre polynomials} 
 $P_{n}(x)$ are examples of Legendre functions, and may be defined via the \emph{Rodrigues formula} 
\begin{equation}\label{201806141227AM1A}
 P_{n}(x) = \frac{1}{2^{n} n!} \frac{d^{n}}{dx^{n}} (x^{2} - 1)^{n}. 
\end{equation}
 The following equivalent definition for $P_{n}(x)$
 will also be used: 
\begin{equation}\label{20190121213PM1A}
 P_{n}(x) = \frac{1}{2^{n}} \sum _{k=0}^n \binom{n}{k}^2 (x-1)^{n-k} (x+1)^k. 
\end{equation}
 The Legendre polynomials form an othogonal family on $(-1, 1)$, with $$ \int_{-1} ^{1} P_{n}(x) P_{m}(x) \, dx = \frac{2}{2n+1} 
 \delta_{m,n}, $$ which gives us the \emph{Fourier--Legendre series} for a suitable function $g$: 
\begin{equation}\label{20190119225PM1A}
 g(x) = \sum_{n=0}^{\infty} \left[ \frac{2n+1}{2} \int_{-1}^{1} g(t) 
 P_{n}(t) \, dt \right] P_{n}(x).
\end{equation}
 Letting shifted Legendre polynomials be denoted as $\tilde{P}_{n}(x) = P_{n}(2x-1)$, polynomials of this form are orthogonal on 
 $[0, 1]$, with $$ \int_{0}^{1} \tilde{P}_{m}(x) \tilde{P}_{n}(x) \, dx = \frac{1}{2n+1} \delta_{m, n}. $$ By analogy with the 
 expansion from \eqref{20190119225PM1A}, for a reasonably well-behaved function $f$ on $(0, 1)$, this function may be 
 expressed in terms of shifted Legendre polynomials by writing $$ f(x) = \sum_{m = 0}^{\infty} c_{m} P_{m}(2x-1). $$ We can 
 determine the scalar coefficient $c_{m}$ 
 in a natural way using the orthogonality of the family of 
 shifted Legendre polynomials. In particular, we see that if we integrate both sides of $$ \tilde{P}_{m}(x) f(x) = \sum_{n = 
 0}^{\infty} c_{n} \tilde{P}_{m}(x) \tilde{P}_{n}(x) $$ over $[0, 1]$, we get $$ c_{m} = (2m+1) \int_{0}^{1} 
 P_{m}(2x-1)\, f(x)\, dx. $$ 
 \emph{Brafman's formula} states that 
\begin{equation}\label{201806091143PM1Ap}
 \sum_{n=0}^{\infty} \frac{(s)_{n} (1-s)_{n}}{(n!)^{2}} P_{n}(x) z^{n} 
 = {}_{2}F_{1}\!\!\left[ 
 \begin{matrix}
 s, 1-s \vspace{1mm}\\ 1 
 \end{matrix} \ \Bigg| \ \alpha \right] 
 {}_{2}F_{1}\!\!\left[ 
 \begin{matrix}
 s, 1-s \vspace{1mm}\\ 1 
 \end{matrix} \ \Bigg| \ \beta \right], 
\end{equation}
 letting $\alpha = \frac{1 - \rho - z}{2}$, $\beta = \frac{1 - \rho + z}{2}$, and $\rho = \sqrt{1 - 2 x z + z^{2}}$ 
 \cite{Brafman1951}.
 The canonical generating function for Legendre polynomials is \cite{MR0058756,MR0064915} 
\begin{equation}\label{20180502851PM1A}
 \frac{1}{\sqrt{1 - 2 x z + z^{2}}} 
 = \sum_{n=0}^{\infty} P_{n}(x) z^{n}. 
\end{equation}
 This gives the following result (see \cite{MR3322254} and \cite{MR0064915}) 
 which we exploit heavily: 
\begin{equation}\label{20190119239PM3A}
 \mathbf{K}\left( \sqrt{k} \right) = 
 \sum_{n \geq 0} \frac{2}{2n+1} P_{n}(2k-1). 
\end{equation} 
 If we make use of the standard moment formula 
\begin{equation}\label{201901211014PM2A}
 \int_0^1 x^i P_{n}(2 x-1) \, dx=\frac{(i!)^2}{(i-n)! (i+n+1)!}
\end{equation} 
 for shifted Legendre polynomials, then it is not difficult to 
 see why \eqref{20190119239PM3A} holds. Namely, 
 from the Maclaurin series for $\text{{\bf K}}$, 
 we have that 
$$\text{{\bf K}}\left(\sqrt{x}\right) P_{n}(2 x-1)
 = \frac{\pi}{2} \sum _{i \in \mathbb{N}_{0}} 
 \left(\frac{1}{16}\right)^i \binom{2 i}{i}^2 x^i P_{n}(2 x-1),$$
 and by rewriting the right-hand side as 
 $$\frac{\pi }{2} \sum _{i \in \mathbb{N}_{0}} \left(\frac{(2 i)!}{i!}\right)^2 \frac{1}{16^i (i-n)! (i+n+1)!}
 = \frac{2 (\sin (\pi n)+1)}{(2 n+1)^2}$$
 using the moments for the family $\{ P_{n}(2x-1) \}_{n \in \mathbb{N}_{0}}$, 
 we obtain the desired result. 

\section{A motivating example}\label{20190119229PM2A}
 Following the integration method from \cite{Campbell2018Ramanujan}, since $$\int_0^1 \frac{x^{4 n} \ln \left(1 - 
 x^2\right)}{\sqrt{1-x^2}} \, dx = -\frac{\sqrt{\pi } \Gamma \left(2 n + \frac{1}{2}\right) \left(H_{2 n} + 
 2 \ln(2)\right)}{2 \Gamma (2 n + 1)},$$ and since $$\sum _{n=0}^{\infty } \frac{(-1)^n x^{4 n} \binom{\frac{1}{2}}{n} 
 \ln \left(1-x^2\right)} {\sqrt{1-x^2}} = \sqrt{x^2+1} \ln \left(1-x^2\right) $$ we have that the series 
\begin{equation*}
\sum _{n=0}^{\infty } \frac{ \binom{2 
 n}{n} \binom{4 n}{2 n} H_{2 n}}{64^{n}(2 n-1)}
\end{equation*}
 may be evaluated in terms of the following ${}_{3}F_{2}(1)$ series, 
 which was conjectured by the first author to be equal to the parbelos constant 
 $ \frac{\sqrt{2} + \ln(1 + \sqrt{2})}{\pi}$ \cite{CampbellDAurizioSondow}: 
\begin{equation*}
 {}_{3}F_{2}\left[ 
 \begin{matrix}
 -\frac{1}{2}, \frac{1}{4}, \frac{3}{4} \\
 \frac{1}{2}, 1 
 \end{matrix} \ \Bigg| \ 1 \right]. 
\end{equation*}
 This conjecture is proved in a variety of different ways in 
 \cite{CampbellDAurizioSondow}, where the ``palindromic'' formula 
\begin{equation}\label{20180315511PM1A}
 {}_{3}F_{2}\!\!\left[ \begin{matrix} 
 \frac{1}{4},\frac{1}{2},\frac{3}{4} \vspace{1mm}\\ 
 1,\frac{3}{2}
 \end{matrix} \ \Bigg| \ 1 \right] = 
 \frac{8}{\pi} \tanh ^{-1} \tan \frac{\pi }{8} 
\end{equation}
\noindent is introduced and is used in one proof 
 through an application of Fourier--Legendre theory 
 that heavily makes use of the complete elliptic integral {\bf K}. 

 Many aspects of our present article are directly inspired by the FL-based proof presented in \cite{CampbellDAurizioSondow}, so it is 
 worthwhile to review it. Adopting notation from \eqref{20190119234PM1A}, we let $g(x) = 
 \frac{1}{ (2 - x)^{3/2}}$; 
 the main integral under investigation is then 
\begin{equation}\label{20190119334PM1A}
 \int_{0}^{1} \text{{\bf K}}\left( \sqrt{x} \right) \frac{1}{\left( 
 2 - x \right)^{3/2}} \, dx. 
\end{equation}
 Letting $(x)_{n} = \frac{\Gamma(x+n)}{\Gamma(x)}$ 
 denote the Pochhammer symbol, by writing the series in 
 \eqref{20180315511PM1A} as 
\begin{align*}
 {}_{3}F_{2}\!\!\left[ \begin{matrix} 
 \frac{1}{4},\frac{1}{2},\frac{3}{4} \vspace{1mm}\\ 
 1,\frac{3}{2}
 \end{matrix} \ \Bigg| \ 1 \right] 
 & = \sum _{n=0}^{\infty } 
 \frac{\left(\frac{1}{4}\right)_n \left(\frac{1}{2}\right)_n 
 \left(\frac{3}{4}\right)_n}{ (1)_n \left(\frac{3}{2}\right)_n n!} \\ 
 & = \sum _{n=0}^{\infty } \frac{\binom{2 n}{n} \binom{4 n}{2 n}}{(2 n+1) 64^n}
\end{align*}
\noindent and by noting that 
 $$\sum _{n=0}^{\infty } \frac{\binom{2 n}{n} \binom{4 n}{2 n} }{64^n} t^{2 n} = 
 \frac{2 \text{{\bf K}}\left(\sqrt{\frac{2 t}{t + 1}}\right)}{\pi \sqrt{t + 1}}, $$
 we thus find that 
\begin{equation}\label{20190119402PM4A}
 {}_{3}F_{2}\!\!\left[ \begin{matrix} 
 \frac{1}{4},\frac{1}{2},\frac{3}{4} \vspace{1mm}\\ 
 1,\frac{3}{2}
 \end{matrix} \ \Bigg| \ 1 \right] 
 = \frac{ 2 \sqrt{2} }{\pi} \int_0^1 \frac{\text{{\bf K}}\left(\sqrt{r}\right)}{(2-r)^{3/2}} \, dr. 
\end{equation}
 Manipulating \eqref{20180502851PM1A} so as to obtain the generating function
 for shifted Legendre polynomials, we find that 
 $$\frac{1}{(2-x)^{\frac{3}{2}}} = \sum_{n \in \mathbb{N}_{0}}
 (2n+1) \sqrt{2} \left( \sqrt{2} - 1 \right)^{2n+1} \tilde{P}_{n}(x),$$
 so that 
 $$ \int_{x \in [0, 1]} \tilde{P}_{m}(x) \frac{1}{(2-x)^{\frac{3}{2}}} \, dx 
 = \sqrt{2} \left( \sqrt{2} - 1 \right)^{2m+1} $$
 for all $m \in \mathbb{N}_{0}$, from the orthogonality relations 
 for shifted FL polynomials. 
 So, from the important expansion in \eqref{20190119239PM3A}, 
 we have that 
 the integral in \eqref{20190119334PM1A} equals 
\begin{align*}
 \int_{0}^{1} \sum_{n \geq 0} \frac{2}{2n+1} \tilde{P}_{n}(x) \frac{1}{\left( 2 - x \right)^{\frac{3}{2}}} \, dx 
 & = \sum_{n \geq 0} \frac{2}{2n+1} \sqrt{2} \left( \sqrt{2} - 1 \right)^{2n+1} \\ 
 & = \sqrt{2} \ln \left(1+\sqrt{2}\right). 
\end{align*}
 This, together with equation \eqref{20190119402PM4A}, gives us the desired evaluation in 
 \eqref{20180315511PM1A}.

 It is natural to consider variants and generalizations of this proof technique, and this serves as something of a basis for our 
 article.

\section{Summary of main results}
 In Section \ref{FLexp}, we show how manipulations related to the fundamental formula \eqref{20190119239PM3A} together 
 with symbolic evaluations of definite integrals of the form given in \eqref{20190119234PM1A} for a function $g(x)$ with a 
 ``reasonable'' or ``manageable'' shifted Fourier--Legendre series may be used to construct new results on hypergeometric 
 series.

 Generalized harmonic functions are of the form 
\begin{equation}\label{20190121152AM2A}
 H_a^{(b)} = \zeta (b) -\zeta (b, a+1),
\end{equation} 
 where $\zeta (b) $ denotes the Riemann zeta function evaluated at $b$, and $$ \zeta\left( b, a+1 \right) = \sum_{i \in 
 \mathbb{N}_{0}} \frac{1}{\left( i + a+1 \right)^{b}} $$ denotes the Hurwitz zeta function with parameters $b$ and $a + 1$. 
 In the case $b = 1$, we often omit the superscript 
 on the left-hand side of \eqref{20190121152AM2A}. We also adopt the standard convention whereby 
 $H_{a}^{(0)} = a$. 
 If we let $c_{n}$, $c_{n}'$, $\ldots$, $c_{n}^{(m)}$ 
 satisfy the hypergeometric condition in \eqref{20190119724PM1A}, 
 then we define a \emph{twisted hypergeometric series}
 to be one of the form 
 $$ \sum_{n=0}^{\infty} \left( c_{n} H_{\alpha n + \beta}^{(\gamma)} 
 + c_{n}' H_{\alpha' n + \beta'}^{(\gamma')} 
 + \cdots + c_{n}^{(m)} H_{\alpha^{(m)} n + \beta^{(m)}}^{(\gamma^{(m)})} \right). $$
 The evaluation of such series using our main technique 
 is of central importance in this paper.

 In Section \ref{20180504539PM1A}, we use this technique to prove that 
\begin{equation}\label{20190121434PM1A}
 \frac{\pi}{4} = \frac{ {}_{3}F_{2}\left[ 
 \begin{matrix} -\eta, \frac{1}{2}, 1 \\ \frac{3}{2}, 2 + \eta \end{matrix} \ \Bigg| \ -1 \right] }{ {}_{3}F_{2}\left[ \begin{matrix} 
 \frac{1}{2}, \frac{1}{2}, 1 + \eta \\ 1, 2 + \eta \end{matrix} \ \Bigg| \ 1 \right]} 
\end{equation}
 for $\eta > -1$, 
 and we obtain the following identity on the moments of the 
 elliptic-type function $\text{{\bf E}}(\sqrt{x})$: 
\begin{align*}
 \int_{0}^{1}\text{{\bf E}}(\sqrt{x})\,x^\eta\,dx 
 & = \frac{\pi}{2(1+\eta)} \cdot 
 {}_{3}F_{2}\left[ \begin{matrix} 
 -\tfrac{1}{2},\tfrac{1}{2},1+\eta \\ 
 1,2+\eta \end{matrix} \ \Bigg| \ 1 \right] \\
 & = \frac{4}{3(1+\eta)}\cdot 
 {}_{3}F_{2}\left[ \begin{matrix} 
 -\tfrac{1}{2},1,-\eta \\ 
 \tfrac{5}{2},2+\eta \end{matrix} \ \Bigg| \ -1 \right]. 
\end{align*}

 In Section \ref{201805011028PM1A}, we prove the equality $$ \sum_{n \geq 0} \binom{4n}{2n} \binom{2n}{n} 
 \frac{ H_{n} - H_{n - 1/2} }{64^{n}} = \frac{\pi}{\sqrt{2}} - \frac{2 \sqrt{2}}{\pi} \ln^{2}\left( \sqrt{2} + 1 \right)$$ and we also 
 provide a closed-form evaluation for the series
\begin{equation*}
\sum _{n=0}^{\infty } \binom{2 n}{n} \binom{4 n}{2 n} \frac{\frac{1}{2}+n 
 H_{n-\frac{1}{2}}-n H_n}{64^{n}}.
\end{equation*}
 We offer a new FL-based proof of the $\frac{1}{\pi}$ formula 
\begin{equation}\label{20180614108PM1A}
\sum_{n=1}^{\infty} \frac{\binom{2n}{n}^{2} H_{n}}{16^{n}(2n-1)} = \frac{8 \ln(2) - 4}{\pi}
\end{equation}
 introduced in 
 \cite{Campbell2018Ramanujan}, along with a new proof of the formula $$ \sum_{m, n \geq 0} \left( 
 \frac{1}{16} \right)^{m} \frac{\binom{2m}{m}^{2} H_{n}}{(n+1)(m+n+2)} = \frac{ 48 + 32(\ln(2) - 2) \ln(2) }{\pi} - \frac{4 \pi } 
{3} $$ given in \cite{Campbell2018Ramanujan}.
 By applying a moment formula that is used to prove 
 the identity in \eqref{20190121434PM1A}, 
 we prove the new double series result 
 $$ \sum _{n, m \in \mathbb{N}_{0}} \frac{ (n+1)! n! H_n }{(2 m+1) (n-m+1)! (m+n+2)!}
 = 12-\frac{\pi ^2}{3}+8 \ln ^2(2)-16 \ln (2).$$
 Inspired in part by an integration method given in \cite{Campbell2018Ramanujan} 
 on the evaluation of series containing factors of the form $H_n^2 + H_n^{(2)}$, we offer a new proof of the formula 
 $$ \sum _{n=0}^{\infty} \left( \frac{1}{16} \right)^{n} \binom{2 n}{n}^2 \frac{ H_n^2 + H_n^{(2)} }{n + 1} = \frac{64 
 \ln^2(2)}{\pi }-\frac{8 \pi }{3} $$ using the machinery of Fourier--Legendre expansions. We also offer a stronger version of this 
 result by showing how our main technique can be used to evaluate $$\sum _{n=0}^{\infty } \left( \frac{1}{16} 
 \right)^{n} \binom{2 n}{n}^2 \frac{ H_n^{(2)}}{n + 1} = \frac{32 G}{\pi }+\frac{2 \pi }{3}-16 \ln (2)$$ which also gives us a 
 symbolic form for its companion $$\sum _{n=0}^{\infty } \left( \frac{1}{16} \right)^{n} \binom{2 n}{n}^2 \frac{ H_n^{2}}{n + 
 1}.$$ We also prove the formula $$ \sum _{n=0}^{\infty } \frac{\binom{2 n}{n}^2 H_{2 n}}{16^n(2 n - 1)} = \frac{6 
 \ln (2)-2}{\pi },$$ which extends our proof of \eqref{20180614108PM1A}. In Section \ref{201805011028PM1A}, we prove 
 the equality 
\begin{equation*}
 \sum _{n=0}^{\infty } \binom{2 n}{n}^2 \frac{H_{n + \frac{1}{4}} - 
 H_{n-\frac{1}{4}}}{16^{n}}
 = \frac{\Gamma^4 \left(\frac{1}{4}\right)}{8 \pi ^2}-\frac{4 G}{\pi}, 
\end{equation*}
and offer a new proof of the equation $$ \sum _{n=1}^{\infty } \frac{ \binom{2 n}{n}^2 H_{2 
 n}}{ 16^{ n} (2 n-1)^2} = \frac{4 G + 6-12 \ln (2)}{\pi }$$ 
 introduced in \cite{CampbellSofo2017}. 
 Using FL theory, we prove the formula 
 $$ \sum _{m, n \geq 0} \frac{\binom{2 m}{m}^2 \binom{2 n}{n}^2 }{ 16^{m + 
 n} (m+n+1) (2 m+3) } = 
 \frac{7 \zeta (3) - 4 G}{\pi ^2}, $$ 
 strongly motivating further explorations on our main techniques. 

 Recall that the polylogarithm function $\text{Li}_{n}(z)$ is defined so that $$\text{Li}_{n}(z) = \sum_{k=1}^{\infty} 
 \frac{z^{k}}{k^{n}}$$ for $|z| \leq 1$, and in the case $n = 2$ we obtain the dilogarithm mapping. In Section 
 \ref{20180609303PM1A} we use our integration methods to explore connections between generalized hypergeometric functions 
 and polylogarithmic functions, providing an evaluation of 
\begin{equation*}
 \sum_{n\geq 1}\frac{\binom{4n}{2n}\binom{2n}{n}}{n\,64^n} 
 =\frac{3}{16}\cdot {}_{4}F_{3}\left[ \begin{matrix}
 1,1,\tfrac{5}{4},\tfrac{7}{4} \\ 
 2,2,2 \end{matrix} \ \Bigg| \ 1 \right] 
\end{equation*}
 in terms of a dilogarithmic expression, with a similar evaluation being given for $$\sum _{n=0}^{\infty } \frac{\binom{2 n}{n}^2 (2 
 n+1)}{16^n (n+1)^4}.$$

 In Section \ref{20180614114PM2A}, we prove a variety of new results on generalized complete elliptic integrals of the form $$ 
 \mathfrak{J}(x) = \int_{0}^{\pi/2}\left(\sqrt{1-x\sin^2\theta}\right)^3\,d\theta.$$ Using the moments of this function, we prove 
 that the identity $$ \frac{15\pi}{32} = \frac{ {}_{3}F_{2}\!\!\left[ \begin{matrix} -\tfrac{3}{2},1,-\eta \vspace{1mm}\\ 
 \tfrac{7}{2},2+\eta \end{matrix} \ \Bigg| \ -1 \right] }{ {}_{3}F_{2}\!\!\left[ \begin{matrix} -\tfrac{3}{2},\tfrac{1}{2},1+\eta 
 \vspace{1mm}\\ 1,2+\eta \end{matrix} \ \Bigg| \ 1 \right]}$$ holds for $\eta > -1$. 
 A generalization $\mathfrak{J}_{m}(x)$ of $ \mathfrak{J}(x)$ 
 is introduced in Section \ref{20180614114PM2A}, and we introduce formulas for evaluating the moments of
 $\mathfrak{J}_{m}(x)$.

\section{Related mathematical literature}
 Some of our main results are the construction of new formulas for $\frac{1}{\pi}$ using Fourier--Legendre expansions. So, it is 
 natural to consider FL-based techniques that have previously been applied to derive new infinite series formulas for 
 $\frac{1}{\pi}$. One of Ramanujan's most famous formulas for $\frac{1}{\pi}$ is $$ \frac{2}{\pi } = \sum _{n=0}^{\infty } 
 \left(-\frac{1}{64}\right)^n (4 n+1) \binom{2 n}{n}^3,$$ which was proved by Bauer in 1859 in \cite{Bauer1859} using a 
 Fourier--Legendre expansion. Levrie applied Bauer's classical FL-based method to functions of the form $(\sqrt{1 - 
 a^2x^2})^{2k-1}$ to derive new infinite sums for $\frac{1}{\pi}$, 
 such as the Ramanujan-like equation \cite{Levrie2010}
 $$\frac{8}{9 \pi } = 
 \sum_{n=0}^{\infty } \left( -\frac{1}{64} \right)^{n} \frac{ (4 n + 1) \binom{2 n}{n}^3}{(n + 1) (n + 2) (2 n - 3) (2 n - 
 1)}. $$ 

 Much of the subject matter in Wan's dissertation \cite{Wan2013} is closely related to 
 some of our main techniques. 
 One of our key methods is the manipulation of generating functions for Legendre polynomials to construct new rational 
 approximations for $\frac{1}{\pi}$, as is the case with \cite{Wan2013}. The section in \cite{Wan2013} on Legendre polynomials 
 and series for $\frac{1}{\pi}$ makes use of Brafman's formula \eqref{201806091143PM1Ap}, proving many series for 
 $\frac{1}{\pi}$ typically involving summands with irrational powers. 

 Wan also explored the use of Legendre polynomials to construct new series for $\frac{1}{\pi}$ in \cite{Wan2014} and new results 
 in this area were also introduced by Chan, Wan, and Zudilin in \cite{ChanWanZudilin2013}. Brafman's formula is also applied in 
 \cite{Wan2014} to produce new results on $\frac{1}{\pi}$ series, whereas the FL-based methods in our present article mainly make 
 use of Fourier--Legendre expansions for elliptic-type expressions such as $\text{{\bf K}}\left( \sqrt{x} \right)$. New series for 
 $\frac{1}{\pi}$ are given in \cite{WanZudilin2012} through a generalization of Bailey's identity for generating functions given by 
 componentwise products of Ap\'{e}ry-type sequences and the sequence of Legendre polynomials. The construction of 
 hypergeometric series identities using expansions in terms of Legendre polynomials has practical applications in mathematical physics 
 \cite{Holdeman1970} and related areas; a variety of binomial sum identities given in terms of generalized hypergeometric functions 
 are proved in \cite{ElCheon2005} through the use of the family $\{ P_{n}(x) : n \in \mathbb{N}_{0} \}$. 

\section{Applications of FL expansions}\label{FLexp}
 From the ordinary generating function for the binomial sequence $$ \left( \binom{4n}{2n} : n \in \mathbb{N}_{0} \right),$$ we 
 have that 
\begin{equation}\label{gf42}
 \sum_{n\geq 0}\binom{4n}{2n}\frac{x^{2n}}{16^n}=\frac{1}{2}\left(\frac{1}{\sqrt{1+x}}+\frac{1}{\sqrt{1-x}}\right).
\end{equation}
 From \eqref{gf42}, together with Wallis' formula
\begin{equation}\label{20190120637PM1A}
 \binom{2n}{n} = \frac{2}{\pi} \int_0^{\frac{\pi }{2}} 2^{2 n} \sin ^{2 n}(t) \, dt, 
\end{equation}
 we may compute the generating function for the sequence $\left(\binom{4n}{2n}\binom{2n}{n}\right)_{n\geq 0}$: 
\begin{align*}
\sum_{n\geq 0}\binom{4n}{2n}\binom{2n}{n}\frac{y^n}{64^n} 
 &= {}_{2}F_{1}\!\!\left[ 
 \begin{matrix}
 \tfrac{1}{4},\tfrac{3}{4} \vspace{1mm}\\ 1 
 \end{matrix} \ \Bigg| \ y \right]\\ 
 &=\frac{1}{\pi}\int_{0}^{\pi/2}\left(\frac{1}{\sqrt{1 - 
 \sqrt{y}\sin\theta}}+\frac{1}{\sqrt{1+\sqrt{y}\sin\theta}}\right)\,d\theta \\
 &= \frac{2}{\pi\sqrt{1+\sqrt{y}}}\,\text{{\bf K}}\left(\sqrt{\frac{2\sqrt{y}}{1+\sqrt{y}}}\right). 
\end{align*}
 The identity \eqref{gf42} alone is powerful enough to lead us to the explicit computation of many hypergeometric 
 $\phantom{}_3 F_2$ functions with quarter-integer parameters. 
 Through the use of \eqref{20190120637PM1A}, 
 it is not difficult to evaluate series of the form 
\begin{equation*}
 \sum_{n\geq 0}\binom{4n}{2n}\binom{2n}{n}\frac{1}{64^n 
 (n+m)^\eta} \qquad m\in\mathbb{N}^+,\eta\in\{1,2\} 
\end{equation*}
 so as to obtain closed-form evaluations for series as in $$ \sum_{n\geq 0}\binom{4n}{2n}\binom{2n}{n} 
 \frac{1}{64^n (n+m_1)^{\eta_1}\cdots(n+m_k)^{\eta_k}} \qquad m_i\in\mathbb{N}^+, \eta_j\in\{1,2\} $$ via partial 
 fraction decomposition, not to mention the case whereby $m\in\tfrac{1}{2}+\mathbb{N}$ and $\eta=1$; in this case we 
 can also apply Wallis' identity to obtain closed-form evaluations. For example, 
 the following equalities hold: 
\begin{align*}
 & \sum_{n\geq 0}\binom{4n}{2n}\binom{2n}{n}\frac{1}{64^n(2n+1)}=\frac{4}{\pi}\ln(1+\sqrt{2}), \\ 
 & \sum_{n\geq 0} 
 \binom{4n}{2n}\binom{2n}{n}\frac{1}{64^n(2n+3)}=\frac{4\sqrt{2}}{15\pi}+\frac{16}{15\pi}\ln(1+\sqrt{2}), \\ 
 & \sum_{n\geq 0}\binom{4n}{2n}\binom{2n}{n}\frac{1}{64^n(2n+5)}=\frac{68\sqrt{2}}{315\pi}+\frac{64}{105\pi}\ln(1+\sqrt{2}). 
\end{align*}

 As Zhou shows in \cite{Zhou2014}, the Legendre functions $P_{-1/2}$ and $P_{-1/4}$ are associated with 
 $\text{{\bf K}}(\sqrt{x})$ and $\frac{1}{\sqrt{1+y}}\text{{\bf K}}\left(\sqrt{\frac{2y}{1+y}}\right)$, which 
 in turn are associated with 
 the weights $\binom{2n}{n}^2$ and $\binom{4n}{2n}\binom{2n}{n}$. 
 From the Rodrigues formula, 
 the computation of the FL expansion of 
 $\frac{1}{\sqrt{1+y}}\text{{\bf K}}\left(\sqrt{\frac{2y}{1+y}}\right)$ turns out to be an exercise in fractional calculus. On the 
 other hand, 
\begin{align*}
 \int_{0}^{1}\frac{y^\eta}{\sqrt{1+y}}\, 
 \text{{\bf K}}\left(\sqrt{\frac{2y}{1+y}}\right)\,dy &= \frac{\pi}{2}\sum_{n\geq 0}\binom{4n}{2n}\binom{2n}{n}\frac{1}{64^n(2n+\eta+1)}
\\&=\frac{\pi}{2(1+\eta)}\cdot 
 {}_{3}F_{2}\left[ 
 \begin{matrix}
 \tfrac{1}{4},\tfrac{3}{4},\tfrac{1+\eta}{2} \\ 
 1,\tfrac{3+\eta}{2}
 \end{matrix} \ \Bigg| \ 1 \right]. 
\end{align*}
 As illustrated in the following example, the main technique in our paper as formulated in the Introduction can also be applied to 
 obtain hypergeometric transformations related to series involving binomial products; we later find that the series in 
 \eqref{20190120155PM1A} arises in our evaluation of a new binomial-harmonic summation. 

\begin{example}
 From the Rodrigues formula \eqref{201806141227AM1A}, we have that: $$\int_0^1 \frac{P_{n}(2 x-1)}{\sqrt{x (1-x)}} \, dx 
 = \frac{ \pi \left((-1)^n+1\right) \binom{n}{\frac{n}{2}}^2}{2^{2 n + 1}}.$$ So, from the fundamental formula 
 \eqref{20190119239PM3A} for expressing $ \mathbf{K}\left( \sqrt{k} \right)$ as a shifted FL expansion, we obtain the following: 
\begin{align*}
 \int_0^1 \frac{\text{{\bf K}}(\sqrt{x})}{\sqrt{(1-x) x}} \, dx 
 & = \int_0^1 \sum_{n \geq 0} \frac{2}{2n+1} \frac{P_{n}(2x - 1)}{\sqrt{(1-x) x}} \, dx \\ 
 & = \sum_{n \geq 0} \frac{2}{2n+1} 
 \frac{ \pi \left((-1)^n+1\right) \binom{n}{\frac{n}{2}}^2}{2^{2 n + 1}} \, dx. 
\end{align*}
 So, the above expressions are equal to 
\begin{equation}\label{20190120155PM1A}
 2 \pi \sum _{n=0}^{\infty } \frac{ \binom{2 n}{n}^2}{16^{n}(4 n + 1)},
\end{equation}
 which \emph{Mathematica} 11 can evaluate directly. 
 On the other hand, by using the Maclaurin series for $\text{{\bf K}}$, 
 we see that 
\begin{equation*}
 \int_0^1 \frac{\text{{\bf K}}(\sqrt{x})}{\sqrt{(1-x) x}} \, dx 
 = \frac{\pi ^2}{2} \sum _{n=0}^{\infty } \left( \frac{1}{64} \right)^{n} \binom{2 n}{n}^3
\end{equation*}
 so that the evaluation as $\frac{\Gamma^4 \left(\frac{1}{4}\right)}{8 \pi }$ 
 follows immediately from the evaluation in terms of $\text{{\bf K}}$ of the generating function 
 for cubed central binomial coefficients.
 Variants of the approach outlined above allow us to compute a wide array of ${}_{p}F_{q}$ series with 
 quarter-integer parameters. 
\end{example}

 From the equalities $$ \mathbf{K}(\sqrt{x}) = \int_{0}^{\pi/2}\frac{d\theta}{\sqrt{1-x \sin^2\theta}} = \frac{\pi}{2}\sum_{n \geq 
 0}\binom{2n}{n}^2\frac{x^n}{16^n}, $$ we see that series of the form 
 $$\sum_{n\geq 0}\frac{\binom{2n}{n}^2}{16^n (n+a)^2} $$ 
 and summations such as 
 $$\sum_{n\geq 0}\frac{\binom{2n}{n}^2 H_{n} }{16^n (n+a)^2} $$ 
 are closely related to 
 the moments of $\text{{\bf K}}\left( \sqrt{x} \right)$ and $\text{{\bf K}}\left( \sqrt{x} \right) \ln(1-x)$, 
 as we later explore. 
 
\subsection{On the moments of elliptic-type integrals}\label{20180504539PM1A}
   The results given in the present subsection are inspired in part by those in \cite{Wan2012}    on the moments of products of 
  complete elliptic integrals.    We remark that the contiguity relations for hypergeometric functions    can be applied to obtain similar 
  results.

\begin{lemma}\label{201805121104PM1A}
  For $m$ such that $\Re(m) > -1$ we have the equation 
\begin{equation}\label{20180608315PM1A}
 \int_0^1 \text{{\bf K}}\left(\sqrt{x}\right) x^m \, dx = 
 2 \sum _{i = 0}^\infty \frac{ (m!)^2}{(2 i + 1) (m-i)! (m+i+1)!}.
\end{equation}
\end{lemma}

\begin{proof}
 Through the use of the Rodrigues formula, it is not difficult to see that the identity $$ \int_0^1 P_{i} (2 x-1) x^m \, dx = 
 \frac{ \Gamma^2 (m + 1)}{\Gamma (m - i + 1) \Gamma (m + i + 2)}$$ holds. So, we have that 
 $$\text{{\bf K}}\left(\sqrt{x}\right) x^m = \sum _{i \in \mathbb{N}_{0}} \frac{ (2 i + 1) (m!)^2 }{ (m-i)! (m + i + 1)!} 
 \text{{\bf K}}\left(\sqrt{x}\right) P_{i}(2 x-1), $$ and the desired result follows by integrating both sides 
 over $[0, 1]$ and applying the fundamental formula \eqref{20190119239PM3A}. 
\end{proof}

\begin{theorem}
 For $\eta > -1$, we have that $$ \frac{\pi}{4} = \frac{ {}_{3}F_{2}\left[ \begin{matrix} -\eta, \frac{1}{2}, 1 \\ \frac{3}{2}, 
 2 + \eta \end{matrix} \ \Bigg| \ -1 \right] }{ {}_{3}F_{2}\left[ \begin{matrix} \frac{1}{2}, \frac{1}{2}, 1 + \eta \\ 1, 2 + \eta 
 \end{matrix} \ \Bigg| \ 1 \right]}. $$ 
\end{theorem}

\begin{proof}
 The right-hand side of formula \eqref{20180608315PM1A} 
 may be evaluated as $$ \frac{2 }{m + 
 1} \cdot {}_{3}F_{2}\!\!\left[ \begin{matrix} \frac{1}{2},1,-m \vspace{1mm}\\ \frac{3}{2}, m 
 + 2 \end{matrix} \ \Bigg| \ -1 \right]. $$ 
 Using the Maclaurin series for $\text{{\bf K}}$
 in the integrand in \eqref{20180608315PM1A}, 
 we see that the 
 definite integral in Lemma \ref{201805121104PM1A}
 must be equal to 
$$\frac{ \pi}{2} \sum _{i \in \mathbb{N}_{0}} \left(\frac{1}{16}\right)^i \frac{\binom{2 i}{i}^2}{i+m+1},$$
 which, in turn, is equal to 
\begin{equation}\label{20190121406PM1A}
 \frac{\pi}{ 2 m+ 2 } \, {}_{3}F_{2}\!\!\left[ \begin{matrix} \frac{1}{2},\frac{1}{2},m+1 \vspace{1mm}\\ 
 1,m+2 \end{matrix} \ \Bigg| \ 1 \right], 
\end{equation}
 giving us the desired result. 
\end{proof}

 Similar methods may be applied to the complete elliptic integral of the second kind, using the FL expansion formula 
\begin{equation}\label{20180609602PM1A}
 \text{{\bf E}}(\sqrt{x}) 
 =4\sum_{n\geq 0}P_n(2x-1)\frac{-1}{(2n-1)(2n+1)(2n+3)},
\end{equation}
\noindent and this leads us to the following result. 

\begin{theorem}\label{22e}
 For $\eta > -1$, we have the equalities
\begin{align*}
 \int_{0}^{1}\text{{\bf E}}(\sqrt{x})\,x^\eta\,dx 
 & = \frac{\pi}{2(1+\eta)} \cdot 
 {}_{3}F_{2}\left[ \begin{matrix} 
 -\tfrac{1}{2},\tfrac{1}{2},1+\eta \\ 
 1,2+\eta \end{matrix} \ \Bigg| \ 1 \right] \\ 
 & = \frac{4}{3(1+\eta)}\cdot 
 {}_{3}F_{2}\left[ \begin{matrix} 
 -\tfrac{1}{2},1,-\eta \\ 
 \tfrac{5}{2},2+\eta \end{matrix} \ \Bigg| \ -1 \right].
\end{align*}
\end{theorem}

\begin{proof}
 The first equality follows from the Taylor series expansion for $\text{{\bf E}}(x)$. Recalling the FL expansion 
 $$ x^\eta = \sum_{n\geq 0}P_n(2x-1)(2n+1)(-1)^n \frac{\Gamma(n-\eta)\Gamma(1+\eta)}{\Gamma(-\eta)\Gamma(n 
 + 2 + \eta)}, $$ we immediately see that $$ \int_{0}^{1} \text{{\bf E}}(\sqrt{x})x^\eta\,dx = 
 \frac{4}{3(1+\eta)}\cdot {}_{3}F_{2}\left[ \begin{matrix} -\tfrac{1}{2},1,-\eta \\ \tfrac{5}{2},2+\eta \end{matrix} \ \Bigg| \ -1 
 \right], $$ as desired. 
\end{proof}

\subsection{New results involving series containing harmonic numbers and central binomial coefficients}\label{201805011028PM1A}
 Inspired in part by the results introduced in 
\cite{Campbell2018New,
Campbell2018Ramanujan,
Campbell2018Series,
CampbellSofo2017}, 
 we evaluate new binomial series for $\frac{1}{\pi}$ involving harmonic numbers, 
 but through the use of Fourier--Legendre theory. 

\begin{theorem}\label{20190210731PM1A}
 The following equality holds: $$ \sum_{n \geq 0} \binom{4n}{2n} \binom{2n}{n} \frac{ H_{n} - H_{n - 1/2} }{64^{n}} = 
 \frac{\pi}{\sqrt{2}} - \frac{2 \sqrt{2}}{\pi} \ln^{2}\left( \sqrt{2} + 1 \right).$$ 
\end{theorem}

\begin{proof}
 From the FL series $$ \frac{1}{\sqrt{2-x}} = \sum_{k \geq 0} 2 (\sqrt{2} - 1)^{2k+1} P_{k}(2x-1) $$ for the elementary 
 function $\frac{1}{\sqrt{2-x}}$, we have that 
 $$ \int_{0}^{1} \frac{ \text{{\bf K}}\left(\sqrt{x}\right) }{\sqrt{2-x}} dx 
 = 4 \sum_{n \geq 0} \frac{ (\sqrt{2} - 1)^{2n+1} }{ (2n+1)^{2}}. $$
 As noted above, the dilogarithm function $\text{Li}_{2}(z)$
 is defined so that 
\begin{equation}\label{20190120910PM2A}
\text{Li}_{2}(z) = \sum_{k \in \mathbb{N}} \frac{z^{k}}{k^{2}}, 
\end{equation}
 whereas the Rogers dilogarithm function $L(x)$ is 
 $$ L(x) = \sum_{n = 1}^{\infty} \frac{x^{n}}{n^{2} } + 
 \frac{1}{2} \ln(x) \ln(1-x). $$ By bisecting the series in \eqref{20190120910PM2A}, we see that $$\sum _{n=0}^{\infty } 
 \frac{x^{2 n+1}}{(2 n+1)^2} = \text{Li}_2(x)-\frac{1}{4}\text{Li}_2\left(x^2\right).$$ Writing $\alpha$ in place of 
 $\sqrt{2} - 1$, we have that 
\begin{equation}\label{20190120951PM2A}
\sum _{n=0}^{\infty } \frac{\alpha^{2 n+1}}{(2 n+1)^2} = 
 \text{Li}_2(\alpha)-\frac{1}{4}\text{Li}_2\left(\alpha^2\right).
\end{equation}
 The evaluation 
 $4 L(\alpha )-L\left(\alpha ^2\right) = \frac{\pi ^2}{4}$
 is known \cite{MR1356515}, which shows that 
$$\text{Li}_2(\alpha )-\frac{\text{Li}_2\left(\alpha ^2\right)}{4}
 = \frac{\pi ^2}{16}-\frac{1}{4} \ln ^2\left(1+\sqrt{2}\right).$$
 So, we have that 
 $$ \int_{0}^{1} \frac{ \text{{\bf K}}\left(\sqrt{x}\right) }{\sqrt{2-x}} dx 
 = \frac{\pi^{2}}{4} - \ln^{2}\left( \sqrt{2} - 1 \right), $$ 
 and by enforcing the substitution $x = \frac{2y}{1 + y}$, we 
 obtain the equality $$ \int_{0}^{1} \text{{\bf K}}\left( \sqrt{\frac{2y}{1+y}} \right) \frac{ \sqrt{2} }{ (1+y)^{3/2} } dy 
 = \frac{\pi}{\sqrt{2}} \sum_{n \geq 0} \binom{4n}{2n} \binom{2n}{n} \frac{1}{64^{n}} \int_{0}^{1} \frac{y^{2n} }{1 + y} 
 dy. $$ 
 Now, we claim that the identity 
 $$ \int_0^1 \frac{y^{2 n}}{1+y} \, dy
 = \frac{H_n-H_{n-\frac{1}{2}}}{2} $$
 holds for $n \in \mathbb{N}_{0}$. 
 Let $ \mathcal{H}_{n}
 = \frac{H_n-H_{n-\frac{1}{2}}}{2} $ 
 for $n \in \mathbb{N}_{0}$, and from 
 \eqref{20190121152AM2A} we see that 
 the recursion
 $$ \mathcal{H}_{n+1} - \mathcal{H}_{n} = 
 -\frac{1}{2 (n+1) (2 n+1)} $$
 holds, with $\mathcal{H}_{0} = \ln(2)$
 as the base case. 
 We see that 
 $ \int_0^1 \frac{y^{2 (0)}}{1+y} \, dy
 = \ln(2), $
 and we also have the equalities 
\begin{align*}
 \int_0^1 \frac{y^{2 (n+1)}}{1+y} \, dy 
 - \int_0^1 \frac{y^{2 n}}{1+y} \, dy 
 & = \int_0^1 \frac{y^{2 (n+1)}-y^{2 n}}{1+y} \, dy \\ 
 & = \int_0^1 (y - 1) y^{2 n} \, dy \\ 
 & = -\frac{1}{2 (n+1) (2 n+1)}.
\end{align*}
\noindent By induction, the desired result follows. 
\end{proof}

\begin{remark}
 \emph{In general, the problem of evaluating a series containing central binomial coefficients and harmonic-type expressions 
 is difficult. However, as explored in \cite{Campbell2018Ramanujan} and \cite{CampbellSofo2017}, there are sometimes 
 ways of obtaining such evaluations through elementary manipulations of generating functions. For example, consider the 
 following variant of the series in Theorem \ref{20190210731PM1A}: 
\begin{equation}\label{201901211237AM1A}
 \sum_{n=0}^{\infty} \binom{2n}{n}^{2} \frac{H_{n} - H_{n - 1/2}}{n+1} \left( 
 \frac{1}{16} \right)^{n}. 
\end{equation}
 In this particular case, if we observe that $$\sum _{n=0}^{\infty } \frac{\binom{2 n}{n} \left(H_n-H_{-\frac{1}{2} + 
 n}\right) x^n}{16^n} = \frac{2 \ln (4-x)}{\sqrt{4-x}} - {}^{\null}_{2}F_{1}^{(0, 0, 1, 0)}\left[ \begin{matrix} \tfrac{1}{2}, 
 1 \\ 1 \end{matrix} \ \Bigg| \ \frac{x}{4} \right] $$ and make use of the identity 
\begin{equation}\label{StandardintegralformulaforCatalannumbers}
 \frac{1}{2 \pi} \int_0^4 x^n 
 \sqrt{\frac{4-x}{x}} \, dx=\frac{\binom{2 n}{n}}{n+1},
\end{equation} 
\noindent it is not difficult to see that the series in \eqref{201901211237AM1A} 
 equals $4-\frac{8}{\pi }$. However, through 
 the use of our main technique, we can also obtain new evaluations that 
 cannot be evaluated in this way, as well as new proofs of known results.}
\end{remark}

 The study of the relationships between definite integrals of the form from \eqref{20190119234PM1A} and generalized 
 hypergeometric functions produces many new results, as further illustrated below. 

\begin{theorem}
 The harmonic-binomial series $$\sum _{n=0}^{\infty } \binom{2 n}{n} \binom{4 n}{2 n} \frac{\frac{1}{2}+n H_{n-\frac{1}{2}}-n 
 H_n}{64^{n}} $$ is equal to $$\frac{3 \pi }{16 \sqrt{2}} + \frac{1}{ 4 \sqrt{2} \pi} \left( 2 + 2 \sqrt{2} \ln \left(1+\sqrt{2}\right) 
 - 3 \ln ^2\left(1+\sqrt{2}\right) \right).$$
\end{theorem} 

\begin{proof} 
 A shifted FL expansion of $\sqrt{2-x}$ is given by 
\begin{equation}\label{20180611552PM1A}
 \sqrt{2-x}=\sum_{n\geq 0}\left(\frac{2}{(1-2n)(2n + 
 3)}+\int_{1}^{2}\frac{dv}{(\sqrt{v}+\sqrt{v-1})^{2n+1}}\right)P_n(2x-1). 
\end{equation}
 Evaluating the integral in \eqref{20180611552PM1A} using the 
 substitution $v\mapsto\cosh^2\theta$, we find 
 that $$ \sqrt{2-x} = \sum_{n\geq 0}\frac{2(\sqrt{2}-1)^{2n+1}(3+\sqrt{2}+2\sqrt{2}\,n)}{(1-2n)(2n+3)}\,P_n(2x-1), $$ and 
 we thus obtain the formula $$ \int_{0}^{1}\text{{\bf K}}(\sqrt{x})\sqrt{2-x}\,dx = 2 + \sqrt{2}\ln(\sqrt{2} - 1) + 12 
 \sum_{n\geq 0}\frac{(\sqrt{2}-1)^{2n+1}}{(1-2n)(2n+1)^2(2n+3)}. $$ Using partial fraction decomposition, 
 we obtain the evaluation 
\begin{equation*}
 \int_{0}^{1}\text{{\bf K}}(\sqrt{x})\sqrt{2-x}\,dx = \frac{8+3 \pi ^2+8 \sqrt{2}\ln(1+\sqrt{2})-12\ln^2(1+\sqrt{2})}{16}. 
\end{equation*}

 If we substitute $x\mapsto\frac{2y}{1+y}$, we see that 
\begin{equation}\label{20190121251AM1A}
 \int_0^1 \text{{\bf K}}\left(\sqrt{x}\right) \sqrt{2-x} \, dx 
 = 2 \sqrt{2} \int_0^1 \text{{\bf K}}\left(\sqrt{\frac{2 y}{1+y}}\right) \left(\frac{1}{1+y}\right)^{5/2} \, dy.
\end{equation}
 Since $$ \text{{\bf K}}\left(\sqrt{\frac{2 y}{1+y}}\right) = \frac{1}{2} \pi \sum _{n=0}^{\infty } \frac{\binom{2 n}{n} 
 \binom{4 n}{2 n} y^{2 n} \sqrt{1+y}}{64^n}, $$ we see that the integral on the right-hand side in 
 \eqref{20190121251AM1A} may be rewritten as $$ \pi\sqrt{2}\sum_{n\geq 0}\binom{4n}{2n}\binom{2n}{n} 
 \frac{1}{64^n}\int_{0}^{1}\frac{y^{2n}}{(1+y)^2}\,dy, $$ and we thus obtain the desired closed form. 
\end{proof}

  The evaluation in Theorem \ref{20180609526PM1A} below    illustrates the power of the integration method from  
  \cite{Campbell2018Ramanujan},    where it is used to give a one-line proof of   \eqref{20180609530PM1A}. Formula   
  \eqref{20180609530PM1A}    is also proved through a complicated manipulation   of the generating function for 
   $\big( \binom{2n}{n} \cdot H_{n} : n \in    \mathbb{N}_{0}\big)$.    Through an application of integration by parts with respect to 
  the known    integral formula $$ \int_{0}^{1} \ln(1-x)\text{{\bf K}}\left(\sqrt{x} \right) \, dx = 8 \ln(2) - 8, $$ we are able to use 
  the    resulting formula    $$ \int_{0}^{1} \frac{1 - \text{{\bf E}}\left( \sqrt{x} \right) }{1 - x} \, dx = 2 - 4 \ln(2) $$   to prove  
  \eqref{20180609530PM1A}, by applying the Maclaurin series for $\text{{\bf E}}\left( \sqrt{x} \right)$. The following proof  
 complements the final proof of \eqref{20180609530PM1A} from \cite{Campbell2018Ramanujan}, and highlights the connections 
 between harmonic-binomial series for $\frac{1}{\pi}$ and Fourier--Legendre theory. 

\begin{theorem}\label{20180609526PM1A}
 The following result holds \cite{Campbell2018Ramanujan}: 
\begin{equation}\label{20180609530PM1A}
\sum_{n=1}^{\infty} \frac{\binom{2n}{n}^{2} H_{n}}{16^{n}(2n-1)} = \frac{8 \ln(2) - 4}{\pi}.
\end{equation}
\end{theorem}

\begin{proof}
 From the Fourier--Legendre series identity 
\begin{equation*}
 \text{{\bf E}}(\sqrt{x}) 
 =4\sum_{n\geq 0}P_n(2x-1)\frac{-1}{(2n-1)(2n+1)(2n+3)},
\end{equation*}
 given in \eqref{20180609602PM1A}, together with the identity $$ \int_{0}^{1} \ln(1 - x) P_{n}\left( 2x - 1 \right) \, dx = 
 - \frac{1}{n(n+1)} $$ we find that 
\begin{equation}\label{20180609705PM1A}
 \int_{0}^{1} \text{{\bf E}}\left(\sqrt{x}\right)\ln(1-x)\,dx 
 = -\frac{4}{3}+\sum_{n\geq 1}\frac{4\left(\frac{1}{n}+\frac{1}{n+1}\right)}{(2n-1)(2n+1)^2 (2n+3)},
\end{equation}
 with the infintie series in \eqref{20180609705PM1A} 
 being evaluated as $ \frac{4}{9}(12\ln 2-11)$. Using the Maclaurin series for 
 $\text{{\bf E}}\left(\sqrt{x}\right)$, we get 
 $$ \sum _{n = 0}^{\infty } \left( \frac{1}{16} 
 \right)^{n} \frac{ \binom{2 n}{n}^2 H_{n + 1}}{(n + 1) (2 n - 1)} = \frac{96 \ln (2)-88}{9 \pi},$$ 
 and the desired result follows using partial fraction decomposition. 
\end{proof}

 As a variation on the main lemma from 
 \cite{Campbell2018Ramanujan}, 
 an integration method in \cite{Campbell2018Ramanujan} 
 allows us to construct new series for $\frac{1}{\pi}$ 
 involving factors of the form $H_{n}^{2} 
 + H_{n}^{(2)}$. For example, the formula
\begin{equation}\label{20180613138PM1A} 
 \sum_{n=1}^{\infty} \frac{ \binom{2n}{n}^{2} \left( H_{n}^{2} + H_{n}^{(2)} \right) }{16^{n}(2n-1)^2} 
 = \frac{64 + 64 \ln^{2}(2) - 96 \ln(2)}{\pi} - \frac{8 \pi}{3}
\end{equation}
 was proved in \cite{Campbell2018Ramanujan} following this method. As noted in \cite{Campbell2018Ramanujan}, 
 using the identity 
\begin{equation*}
 \int_{0}^{1} n x^{n-1} \ln^{2}(1-x) \, dx = (H_{n})^{2} + H_{n}^{(2)}, 
\end{equation*}
 we find that the formula in 
 \eqref{20180613138PM1A} is equivalent to the evaluation 
\begin{equation}\label{20180613142PM1A}
 \int_{0}^{1} \text{{\bf K}}\left( \sqrt{x} \right) \ln^{2}\left( 1 - x \right) \, dx 
 = 48 - \frac{4 \pi^{2}}{3} + 32 (\ln(2) - 2) \ln(2). 
\end{equation}
 As also noted in \cite{Campbell2018Ramanujan}, from \eqref{20180613138PM1A} and \eqref{20180613142PM1A}, 
 together with the identity 
\begin{equation}\label{20190121430PM2A}
 \sum_{n=0}^{\infty} \frac{x^{n+1} H_{n}}{n+1} = \frac{1}{2} \ln^{2}(1-x) 
\end{equation} 
we obtain the 
 double series formula given in 
 Theorem \ref{20180613240PM1A} below, 
 which we prove using the machinery of FL expansions, 
 thereby furthering the main ``thesis'' of our article concerning 
 the connections among the evaluation of twisted and 
 non-twisted hypergeometric series, that of integrals involving elliptic-type integrals, and Fourier--Legendre theory. 

\begin{lemma}\label{20180613146PM1A}
 For $n \in \mathbb{N}$ we have the identity 
\begin{equation}\label{20190121215PM2A}
\int_0^1 \ln ^2(1-x) P_{n}(2 x-1) \, dx
 = \frac{4 n + 2}{n^2 (n + 1)^2}+\frac{4 H_{n - 1}}{n (n + 1)}. 
\end{equation}
\end{lemma}

\begin{proof}
 From \eqref{20190121213PM1A}, we see that
 the equality in \eqref{20190121215PM2A}
 is equivalent to 
\begin{equation}\label{20190121229PM1A}
 \sum _{k=0}^n (-1)^{n - k} \binom{n}{k} \left(\left(H_{n + 1}-H_{n - k}\right){}^2+H_{n + 1}^{(2)}-H_{n 
 - k}^{(2)}\right) 
\end{equation}
 being equal to $\frac{4 n + 2 }{n^2 (n + 1)}+\frac{4 H_{n - 1}}{n}$. 
 Expanding the summand in \eqref{20190121229PM1A} and simplifying using the binomial theorem, 
 we obtain the sum 
\begin{align*}
 & \sum _{k=0}^n (-1)^{n - k + 1} \left(2 \binom{n}{k} H_{n + 1} H_{n - k} - 
 \binom{n}{k} H_{n - k}^2 + \binom{n}{k} H_{n - k}^{(2)} \right). 
\end{align*}
 The following basic binomial-harmonic identities 
 are known to hold: 
\begin{align*}
 \sum _{k=0}^n (-1)^k \binom{n}{k} H_{n-k} & = \frac{(-1)^{n + 1}}{n}, \\ 
 \sum _{k=0}^n (-1)^{n-k} \binom{n}{k} H_{n-k}^2 & = \frac{H_n}{n}-\frac{2}{n^2}, \\ 
 \sum _{k=0}^n (-1)^{n-k+1} \binom{n}{k} H_{n-k}^{(2)} & = \frac{1+n H_{n-1}}{n^2}. 
\end{align*}
 From them, together with the above expansion of 
 \eqref{20190121229PM1A}, we see that the desired result holds.
\end{proof}

\begin{theorem}\label{20180613240PM1A}
 The following equality holds \cite{Campbell2018Ramanujan}: 
 $$ \sum_{m, n \geq 0} \left( \frac{1}{16} \right)^{m} 
 \frac{\binom{2m}{m}^{2} H_{n}}{(n+1)(m+n+2)} = \frac{ 48 
 + 32(\ln(2) - 2) \ln(2) }{\pi} - \frac{4 \pi }{3}. $$
\end{theorem}

\begin{proof}
 By Lemma \ref{20180613146PM1A}, we have that 
\begin{equation}\label{20180613224PM1A}
 \int_{0}^{1} \text{{\bf K}}\left(\sqrt{x}\right)\ln^2(1-x)\,dx 
 = 1 + \sum_{n\geq 1}\frac{1}{n^2(n+1)^2} + 2 \sum_{n\geq 1}\frac{H_{n-1}}{n(n+1)(2n+1)}, 
\end{equation}
 and we thus find that an application of partial fraction decomposition shows that \eqref{20180613142PM1A} holds, 
 giving us a new proof of it. 
 Rewriting the integral in \eqref{20180613224PM1A} as 
 $$\int_0^1 \text{{\bf K}}\left(\sqrt{x}\right) \sum _{n=0}^{\infty } \frac{2 x^{n + 1} H_n}{n + 1} \, dx,$$
 and from the hypergeometric evaluation in \eqref{20190121406PM1A} 
 for the moments indicated in Lemma \ref{201805121104PM1A}, 
 we obtain the desired result.
\end{proof}

 On the other hand, if we use the alternative moment formula in Lemma \ref{201805121104PM1A}, by mimicking the above 
 proof, we obtain the following evaluation.

\begin{theorem}
 The following equality holds: 
 $$ 12-\frac{\pi ^2}{3}+8 \ln ^2(2)-16 \ln (2) 
 = \sum _{n, m \in \mathbb{N}_{0}} \frac{ (n+1)! n! H_n }{(2 m+1) (n-m+1)! (m+n+2)!}. $$
\end{theorem}

\begin{proof}
 This follows almost immediately from
 Lemma \ref{201805121104PM1A} 
 and formulas \eqref{20180613142PM1A} and \eqref{20190121430PM2A},
 in the manner suggested above.
\end{proof}

 The following evaluation may also be obtained using the latter integration method from \cite{Campbell2018Ramanujan}, but the 
 following FL proof illustrates the connections between this evaluation and recent results from \cite{Campbell2018New}.

\begin{theorem}\label{20181018618PM1A}
 The following equality holds: 
 $$ \sum _{n=0}^{\infty } \left( \frac{1}{16} \right)^{n} \binom{2 n}{n}^2 
 \frac{ H_n^2 + H_n^{(2)} }{n + 1} = 
 \frac{64 \ln^2(2)}{\pi }-\frac{8 \pi }{3}. $$ 
\end{theorem}

\begin{proof}
 Expanding the left-hand side of \eqref{20180613224PM1A}, we find that $$ {\frac{4}{3} \left(36-\pi ^2-48 \ln(2) + 
 24\ln^2(2)\right)} = \frac{\pi}{2}\sum_{n\geq 0}\binom{2n}{n}^2\frac{H_{n+1}^2+H_{n+1}^{(2)}}{16^n (n+1)}.$$ 
 By expanding the numerator in 
 the summand, and by evaluating the generalized hypergeometric expression 
 $$ {}_{4}F_{3}\!\!\left[ \begin{matrix} \frac{1}{2},\frac{1}{2},1,1 \vspace{1mm}\\ 2,2,2 \end{matrix} \ \Bigg| \ 1 \right] 
 = \frac{16 (-2 G + 3+\pi (\ln (2)-1))}{\pi } $$ using 
 the standard integral formula \eqref{StandardintegralformulaforCatalannumbers} 
 for the Catalan number $C_{n} = \frac{\binom{2n}{n}}{n+1}$, 
 and by making use 
 of the recently-discovered formula $$ \sum_{n=0}^{\infty} \frac{ \binom{2n}{n}^{2} H_{n} }{16^{n} (n+1)^{2}} = 16 + 
 \frac{32 G - 64 \ln(2)}{\pi} - 16 \ln(2) $$ 
 introduced in \cite{Campbell2018New}, we obtain the desired result. 
\end{proof}

 Through the use of 
 our main technique, as applied to the shifted FL expansion for $\ln(1-x) 
 \ln(x)$, we obtain the following stronger version of 
 Theorem \ref{20181018618PM1A}, which also allows us to evalate $$ \sum_{n \in 
 \mathbb{N}_{0}} \left( \frac{1}{16} \right)^{n} \binom{2n}{n}^{2} \frac{H_{n}^{2} }{n+1} $$ in closed form. The 
 proof depends in a non-trivial way on a recent evaluation provided in \cite{Campbell2018New}. 

\begin{theorem}
 The following equality holds: 
 $$\sum _{n=0}^{\infty } \left( \frac{1}{16} \right)^{n} \binom{2 n}{n}^2 \frac{ H_n^{(2)}}{n + 1} = \frac{32 
 G}{\pi }+\frac{2 \pi }{3}-16 \ln (2). $$ 
\end{theorem}

\begin{proof}
 We find that the equality 
 $$-\frac{\left((-1)^m+1\right) (2 m+1)}{m^2 (m+1)^2} 
 = (2 m+1) \int_0^1 P(m,2 x-1) \ln(1-x) \ln(x) \, dx$$ 
 holds for $m \in \mathbb{N}$.
 By expressing the integral 
 in \eqref{20190119234PM1A}
 in the two different ways our main strategy is based upon, 
 letting $g(x) = \ln(1-x) \ln(x)$ in \eqref{20190119234PM1A}, 
 we obtain the equality 
 $$\sum _{n=0}^{\infty } \frac{ \binom{2 n}{n}^2 \left(1+(1+n) H_{1+n}-(1+n)^2 \psi ^{(1)}(1+n)\right)}{16^n (n + 1)^3}
 = \frac{96-64 \ln (2)}{\pi }-16,$$ 
 and from the recent evaluation \cite{Campbell2018New} 
 $$ \sum_{n \in \mathbb{N}} \frac{ \binom{2n}{n}^{2} H_{n} }{ 16^{n} (n+1)^{2} } 
 = 16 + \frac{32 G - 64 \ln(2)}{\pi} - 16 \ln(2),$$
 we obtain the desired result.
\end{proof}

 The integration result in the following lemma 
 is especially useful in the construction of harmonic-binomial series. 
 We may obtain similar results for integrals involving 
 expressions such as $\ln(1+\sqrt{1-x})$, $ \ln(1+\sqrt{x})$, 
 $ \ln(1-\sqrt{1-x})$, etc.

\begin{lemma}\label{201901211109PM1A}
 The following identity holds for $n \in \mathbb{N}$: 
\begin{equation}\label{201901211040PM1A}
 \int_0^1 \ln \left(1-\sqrt{x}\right) P_{n}(2 x-1) \, dx=\frac{(-1)^n - 4 n - 2}{2 n (n + 1) (2 n + 1)}. 
\end{equation}
\end{lemma}

\begin{proof}
 It seems natural to make use of the moment formula 
 from \eqref{201901211014PM2A}, 
 together with a Maclaurin expansion 
 of $\ln(1 - \sqrt{x})$. 
 Indeed, since 
 $$\sum _{i=0}^{\infty } -\frac{x^{\frac{1+i}{2}} P_{n}(2 x-1)}{1+i}=\ln \left(1-\sqrt{x}\right) P_{n}(2 x-1),$$ 
 from \eqref{201901211014PM2A} we see that 
 $$\int_0^1 \ln \left(1-\sqrt{x}\right) P_{n}(2 x-1) \, dx 
 = - \sum _{i=0}^{\infty } \frac{\left( \left( \frac{i + 1}{2}\right)!\right)^2}{(i + 1) 
 \left(\frac{i + 1}{2}-n\right)! \left(\frac{ i + 1} {2}+n+1\right)! }.$$
 Rewrite the series as 
 $$-\sum _{i=0}^{\infty } \frac{\Gamma^2 \left(\frac{i + 3}{2}\right)}{(i + 1) 
 \Gamma \left(\frac{i}{2} - n + \frac{3}{2} \right) \Gamma \left( \frac{i}{2} + n + \frac{5}{2}\right)}.$$
 Bisect this resulting series so as to produce the symbolic form 
 $$\frac{-2 n+\cos (\pi n)-1}{4 n^3+6 n^2+2 n}-\frac{\pi n (n+1)-\sin (\pi n)}{2 \pi n^2 (n+1)^2},$$
 which is the same as the right-hand side of \eqref{201901211040PM1A}
 for the desired integer parameters. 
\end{proof}

 The formula in Lemma \ref{201901211109PM1A} 
 can be used in conjunction with the identity $\int_{0}^{1} x^n\ln(1-\sqrt{x})\,dx = 
 -\frac{H_{2n+2}}{n+1}$ to produce some new results. For example, since 
\begin{equation}\label{201805011042PM1A}
 \int_{0}^{1} \text{{\bf K}}\left( \sqrt{x} \right)\ln\left(1-\sqrt{x}\right)\,dx 
 = -\frac{\pi}{2}\sum_{n\geq 0}\frac{\binom{2n}{n}^2 H_{2n+2}}{16^n (n+1)},
\end{equation}
\noindent by rewriting the integral in \eqref{201805011042PM1A} 
 as $$-3+\sum_{n\geq 1}\left[\frac{(-1)^n}{2}\left(\frac{1}{n} - 
 \frac{1}{n+1}\right)-\left(\frac{1}{n}+\frac{1}{n+1}\right)\right]\frac{2}{(2n+1)^2},$$ through an application of partial fraction 
 decomposition we obtain a new proof of the known series formula 
 $$ \sum _{n=0}^{\infty } 
 \frac{\binom{2 n}{n}^2 H_{2 n}}{16^n (n+1)} 
 = 2+\frac{4-12 \ln (2)}{\pi } $$ 
 which \emph{Mathematica} 11 is able to verify. However, if we apply the same 
 approach with respect to the definite integral obtained by replacing $\text{{\bf K}}$ in the integrand in 
 \eqref{201805011042PM1A} with the complete elliptic integral of the second kind, we obtain an evaluation of a new 
 binomial-harmonic series. 

\begin{theorem}\label{20180614331AM1A}
 We have the evaluation 
 $$ \sum _{n=0}^{\infty } \frac{\binom{2 n}{n}^2 H_{2 n}}{16^n(2 n - 1)} 
 = \frac{6 \ln (2)-2}{\pi }.$$ 
\end{theorem} 

\begin{proof}
 Rewrite the integral 
\begin{equation}\label{201901211049PM1A}
 \int_{0}^{1} \text{{\bf E}}\left( \sqrt{x} \right) \cdot \ln(1-\sqrt{x}) \, dx 
\end{equation}
 using the Maclaurin series for $\text{{\bf E}}$, 
 so that the integral is equal to 
 $$\frac{\pi}{2} \sum _{i=0}^{\infty } 
 \left(\frac{1}{16}\right)^i \frac{ \binom{2 i}{i}^2 H_{2 i+2} }{(2 i-1) (i + 1)}.$$
 The theorem then follows naturally 
 from Lemma \ref{201901211109PM1A}, 
 as the integral in \eqref{201901211049PM1A}
 may also be written as 
 $$ -2+\sum_{n\geq 1}\left[\frac{(-1)^n}{2}\left(\frac{1}{n}-\frac{1}{n+1}\right)-\left(\frac{1}{n} + 
 \frac{1}{n+1}\right)\right]\frac{4}{(1-2n)(2n+1)^2(2n+3)}. $$ 
 Using partial fraction decomposition, we obtain the desired result. 
\end{proof}

\begin{theorem}
 The following equality holds: 
 $$ \sum _{n=0}^{\infty } \binom{2 n}{n}^2 
 \frac{H_{n + \frac{1}{4}}-H_{n-\frac{1}{4}}}{16^{n} } 
 = \frac{\Gamma^4 \left(\frac{1}{4}\right)}{8 \pi ^2}-\frac{4 G}{\pi}. $$ 
\end{theorem}

\begin{proof} 
 By recalling two standard (equivalent) definitions $$ G = \frac{1}{2}\int_{0}^{\pi/2}\frac{\theta\,d\theta}{\sin\theta}, 
 \qquad G=\sum_{n\geq 0}\frac{(-1)^n}{(2n+1)^2} $$ of Catalan's constant, we may prove (by exploiting the 
 FL expansions of $\text{{\bf K}}(\sqrt{x})$ and $\frac{1}{\sqrt{x}}$, together with the substitution $\theta \mapsto 
 \arcsin\sqrt{t}$) the following identity: $$ 4G = \int_{0}^{1}\frac{\arcsin\sqrt{t}}{\sqrt{t}} \cdot \frac{1}{\sqrt{t (1 - t)}}\,dt 
 = \int_{0}^{1}\frac{\text{{\bf K}}(\sqrt{x})}{\sqrt{x}}\,dx. $$ From the FL expansion $$ \frac{\arcsin\sqrt{x}}{\sqrt{x}} = 
 \sum_{n\geq 0}P_n(2x - 1)\left[\frac{2}{2n+1}-\int_{0}^{1} \frac{4x^{2n + 2}}{1 + x^2}\,dx\right] $$ 
 together with the evaluation for \eqref{20190120155PM1A}, 
 we find that $$ 
 4G = \underbrace{2\pi\sum_{n\geq 0}\frac{\binom{2n}{n}^2}{16^n(4n+1)}}_{\frac{1}{8\pi}\Gamma^2 
 \left(\frac{1}{4}\right)}-4\pi\int_{0}^{1}\sum_{n\geq 0}\frac{\binom{2n}{n}^2}{16^n}x^{4n+2}\frac{dx}{1+x^2}$$ and that 
\begin{align*} G &= \int_{0}^{1}\frac{1-x^2}{1 + 
 x^2}\, \text{{\bf K}}(x^2)\,dx 
 = \frac{1}{4}\int_{0}^{1}\frac{1-\sqrt{x}}{1+\sqrt{x}}\,x^{-3/4} \text{{\bf K}}(\sqrt{x})\,dx\\ 
&= \frac{\pi}{4}\sum_{n\geq 0}\frac{\binom{2n}{n}^2}{16^n}\left(\frac{2}{4n+1}+H_{n-1/4}-H_{n+1/4}\right)
\end{align*}
\noindent as desired.
\end{proof}

 Given the FL expansion of $x\,\text{{\bf K}}(\sqrt{x})$, and introducing the notation 
 $$W(m)=\frac{1}{2m-1}+\frac{2}{2m+1}+\frac{1}{2m+3}+\frac{1}{(2m-1)^2}-\frac{1}{(2m+3)^2}$$ 
for the sake of brevity, we may compute the following integral:
\begin{align*}\int_{0}^{1}\frac{x\, \text{{\bf K}}(\sqrt{x})}{\sqrt{1-x}}\,dx 
 &= \sum_{m\geq 0}\frac{W(m)}{2m+1}\\
&= \frac{3\pi^2}{8} = 2\pi\sum_{n\geq 0}\frac{\binom{2n}{n}(n+1)}{4^n(2n+1)(2n+3)}\\
&= \frac{\pi}{2}\sum_{n\geq 0}\frac{\binom{2n}{n}}{4^n(2n+1)}+\frac{\pi}{2}\sum_{n\geq 0}\frac{\binom{2n}{n}}{4^n(2n 
 + 3)}. 
\end{align*}

 Using the sequence $(W(m) : m \in \mathbb{N}_{0})$, we offer a new proof of the following evaluation 
 introduced in \cite{CampbellSofo2017}. 

\begin{theorem}\label{20180614127AM1A}
 The following equality holds \cite{CampbellSofo2017}: 
 $$ \sum _{n=1}^{\infty } \frac{ \binom{2 n}{n}^2 H_{2 n}}{ 16^{ n} (2 n-1)^2} 
 = \frac{4 G + 6-12 \ln (2)}{\pi}. $$ 
\end{theorem}

\begin{proof}
 Letting $W$ be as given above, we find that 
\begin{align*}
 & \int_{0}^{1}x\, \text{{\bf K}}(\sqrt{x})\ln(1-\sqrt{x})\,dx \\ 
 & = -\frac{5}{3}+\frac{1}{2}\sum_{m\geq 1}\frac{W(m)}{2m+1}\left[\frac{(-1)^m}{2}\left(\frac{1}{m}-\frac{1}{m+1}\right)-\left(\frac{1}{m}+\frac{1}{m+1}\right)\right].
\end{align*}
 Evaluate the right-hand side as $$ \frac{10}{3}\ln(2)-\frac{7}{3}G-\frac{151}{54}. $$ 
 Evaluate the above integral as 
 $$ -\frac{\pi}{2}\sum_{n\geq 0}\frac{\binom{2n}{n}^2 H_{2n+4}}{16^n(n+2)}.$$ 
 It is not difficult to see that it follows that 
 $$\sum _{n = 1}^{\infty } \left( \frac{1}{16} \right)^{n} \frac{ \binom{2 n}{n}^2 H_{2 n}}{n+2} = \frac{92+24 \pi 
 -180 \ln (2)}{27 \pi }.$$ Through a re-indexing argument, it is seen that $$ \frac{12 G+8-6 \ln (2)}{9 \pi } = \frac{1}{3} 
 \sum _{n=1}^{\infty } \frac{ \binom{2 n}{n}^2 H_{2 n}}{16^{ n}(2 n-1)^2}+\frac{5}{9} \sum _{n=1}^{\infty } \frac{ \binom{2 
 n}{n}^2 H_{2 n}}{ 16^{n} (2 n-1)} $$ and the desired result follows from Theorem \ref{20180614331AM1A}. 
\end{proof}

\subsection{Double hypergeometric series}
 We prove a \emph{transformation formula} related to the main results from the preceding 
 sections, and we apply this formula to arrive at a new evaluation for a double hypergeometric series. 

\begin{lemma}\label{20180613619PM1A}
 If $g\in L^1(0,1)$, then the definite integral $$\int_{0}^{1} \text{{\bf E}}(\sqrt{x}) g(x)\,dx$$ is equal to $$ \int_{0}^{1} 
 \left[\frac{\pi}{2}-\left(\frac{\pi}{2}-1\right)\sqrt{x}\right]g(x)\,dx + \frac{1}{3}\iint_{(0,1)^2} 
 \left(\text{{\bf K}}(\sqrt{x})-\frac{\pi}{2}\right)\,g(xz^{2/3})\,dz\,dx. $$
\end{lemma}

\begin{proof}
 Write $ \tilde{g}(x) = \int_{0}^{x} z^2 g(z^2)\,dz$. We note that $$ \frac{\tilde{g}(x)}{x\sqrt{x}} = \frac{1}{x\sqrt{x}} 
 \int_{0}^{\sqrt{x}} z^2 g(z^2)\,dz = \int_{0}^{1} z^2 g(xz^2)\,dz.$$ Using the Maclaurin series expansion for the complete 
 elliptic integral of the second kind, we find that $$ \int_{0}^{1}\text{{\bf E}}(\sqrt{x}) g(x)\,dx = \frac{\pi}{2} 
 \int_{0}^{1}g(x)\,dx - \pi\int_{0}^{1}\sum_{n\geq 1}\frac{\binom{2n}{n}^2 x^{2n-1}}{16^n(2n-1)} x^{2} g(x^2)\,dx. $$ Using 
 integration by parts, we may rewrite the right-hand side as $$ \frac{\pi}{2}\int_{0}^{1}g(x)\,dx - 
 \pi\left[\sum_{n\geq 1}\frac{\binom{2n}{n}^2 x^{2n-1}}{16^n(2n-1)}\tilde{g}(x)\right]_{0}^{1} +\pi\int_{0}^{1}\sum_{n \geq 
 1}\frac{\binom{2n}{n}^2 x^{2n}}{16^n}\cdot\frac{\tilde{g}(x)}{x^2}\,dx.$$ Rewriting this as 
\begin{equation}\label{20180613558PM1A}
 \frac{\pi}{2}\int_{0}^{1}g(x)\,dx -(\pi-2)\int_{0}^{1}z^2 g(z^2)\,dz 
 + \int_{0}^{1}\left(\text{{\bf K}}(\sqrt{x})-\frac{\pi}{2}\right)\frac{\tilde{g}(\sqrt{x})}{x\sqrt{x}}\,dx,
\end{equation}
 we see that \eqref{20180613558PM1A} is equal to 
\begin{align}
 & \int_{0}^{1}\left[\frac{\pi}{2}-\left(\frac{\pi}{2}-1\right)\sqrt{x}\right]g(x)\,dx+ \int_{0}^{1}\left(\text{{\bf K}}(\sqrt{x})-\frac{\pi}{2}\right)\int_{0}^{1} z^2 g(xz^2)\,dz\,dx\nonumber\\
 &= \int_{0}^{1}\left[\frac{\pi}{2}-\left(\frac{\pi}{2}-1\right)\sqrt{x}\right]g(x)\,dx+ \frac{1}{2}\iint_{(0,1)^2}\left(\text{{\bf K}}(\sqrt{x})-\frac{\pi}{2}\right)\sqrt{z} g(xz)\,dz\,dx\nonumber
\end{align}
\noindent and this gives us the desired result. 
\end{proof}

\begin{theorem}
 The following formula holds: 
 $$ \sum _{m, n \geq 0} \frac{\binom{2 m}{m}^2 \binom{2 n}{n}^2 }{ 16^{m + n} 
 (m+n+1) (2 m+3) } = 
 \frac{7 \zeta (3) - 4 G}{\pi ^2}. $$ 
\end{theorem}

\begin{proof}
 We begin with the following:
\begin{align*}
 S & = \sum_{m,n\geq 0}\frac{\binom{2n}{n}^2 \binom{2m}{m}^2}{16^{n+m}(n+m+1)(2m+3)} \\ 
 & = \iint_{(0,1)^2}\sum_{m,n\geq 0}\frac{\binom{2n}{n}^2 
 \binom{2m}{m}^2 x^{n+m} y^{2m+2}}{16^{n+m}}\,dx\,dy.
\end{align*}
 Due to the Maclaurin series of $K(x)$ the above expressions are equal to:
 $$ \frac{4}{\pi^2}\iint_{(0,1)^2} y^2 \text{{\bf K}}(\sqrt{x}) \text{{\bf K}}(\sqrt{x} y)\,dx \,dy. $$
By integrating with respect to $dy$ first we have
$$ S = \frac{2}{3\pi}\int_{0}^{1} 
 \text{{\bf K}}(\sqrt{x}) \cdot 
 \underbrace{ {}_{3}F_{2}\!\!\left[ \begin{matrix} \frac{1}{2}, \frac{1}{2}, \frac{3}{2} \vspace{1mm}\\ 
 1 , \frac{5}{2} \end{matrix} \ \Bigg| \ x \right] }_{g(x)}\,dx, $$
 with $$ g(x) = \sum_{n\geq 0}\left[\frac{1}{4^n}\binom{2n}{n}\right]^2 \frac{3x^n}{2n+3}. $$
 Lemma \ref{20180613619PM1A} gives us that 
\begin{align*}
 & \iint_{(0,1)^2}\left(\text{{\bf K}}(\sqrt{x})-\frac{\pi}{2}\right)y^2 g(xy^2)\,dx\,dy \\ 
 & = \int_{0}^{1} \text{{\bf E}}(\sqrt{x})g(x)\,dx - 
 \int_{0}^{1}\sqrt{x} g(x)\,dx - \frac{\pi}{2}\int_{0}^{1}(1-\sqrt{x})g(x)\,dx. 
\end{align*}
 By considering the instance $g(x)=\text{{\bf K}}(\sqrt{x})$ we have 
\begin{align*}
 & \iint_{(0,1)^2}\left(\text{{\bf K}}(\sqrt{x}) - 
 \frac{\pi}{2}\right)y^2 \text{{\bf K}}(\sqrt{x} y)\,dx\,dy \\ 
 & = \frac{2+7\zeta(3)}{4}-\frac{1+2G}{2}-\frac{\pi(3-2G)}{4}. 
\end{align*}
 The above evaluation is obtained 
 through the FL-expansions of $\text{{\bf K}}(\sqrt{x})$, $\text{{\bf E}}(\sqrt{x})$, 
 and $\sqrt{x}$. 

 In order to complete the evaluation of $S$ it is enough to compute 
\begin{align*}
 \iint_{(0,1)^2} y^2 \text{{\bf K}}(\sqrt{x} y ) \,dx\,dy 
 & = \frac{\pi}{6}\int_{0}^{1}g(x)\,dx \\ 
 & = \frac{\pi}{2}\sum_{n\geq 0}\left[\frac{1}{4^n}\binom{2n}{n}\right]^2 \frac{3}{(2n+3)(n+1)}, 
\end{align*}
 but this is straightforward by partial fraction decomposition, since 
\begin{align*}
 \frac{\pi}{2}\sum_{n\geq 0}\left[\frac{1}{4^n}\binom{2n}{n}\right]^2 \frac{1}{2n+3} 
 & = \int_{0}^{1} x^2 \text{{\bf K}}(x)\,dx \\ 
 & = \frac{1}{2}\int_{0}^{1}\sqrt{x} \text{{\bf K}}(\sqrt{x})\,dx \\ 
 & = \frac{1+2G}{4},
\end{align*}
 and since 
 $$\frac{\pi}{2}\sum_{n\geq 0}\left[\frac{1}{4^n}\binom{2n}{n}\right]^2 \frac{1}{n+1} 
 = \int_{0}^{1} \text{{\bf K}}(\sqrt{x})\,dx = 2. \qedhere$$ 
\end{proof}

 We have that $$ \int_{-1}^{1}x P_N(x)P_L(x)\,dx=\left\{\begin{array}{rcl}\frac{2L+2}{(2L+1)(2L+3)}&\text{if}& N=L+1\\ 
 \frac{2L}{(2L-1)(2L+1)}&\text{if}&N=L-1,\end{array}\right. $$ and hence $$ \int_{0}^{1}x P_N(2x - 1)P_L(2 x - 1)\,dx = 
 \left\{\begin{array}{rcl}\frac{2L+2}{4(2L+1)(2L+3)}&\text{if}& N=L+1\\ \frac{1}{2(2L+1)}&\text{if}& N=L\\ \frac{2L}{4(2L - 
 1) (2L+1)}&\text{if}&N=L-1\end{array}\right. $$ so the previously computed FL expansions allow us an explicit evaluation of 
 any integral of the form $\int_{0}^{1} x \text{{\bf K}}(\sqrt{x}) g(x)\,dx$ or $\int_{0}^{1} x \text{{\bf E}}(\sqrt{x}) g(x)\,dx$, 
with $g(x)$ being a function with a previously computed FL expansion. 

\subsection{Connections between generalized hypergeometric functions and polylogarithms}\label{20180609303PM1A}
 The FL-based techniques explored in our article are powerful enough to find symbolic expressions for some $\phantom{}_4 F_3$ 
 functions with quarter-integer parameters, as shown in the proof of the following theorem. Our manipulations of Fourier--Legendre 
 expansions often allow us to express ${}_{p}F_{q}$ 
 series in terms of the dilogarithm, 
 expanding upon results given in \cite{Daurizio2017}. 

\begin{theorem}
 The hypergeometric series 
\begin{equation}\label{20180609433PM1A}
 \sum_{n\geq 1}\frac{\binom{4n}{2n}\binom{2n}{n}}{n\,64^n} 
 =\frac{3}{16}\cdot {}_{4}F_{3}\left[ \begin{matrix}
 1,1,\tfrac{5}{4},\tfrac{7}{4} \\ 
 2,2,2 \end{matrix} \ \Bigg| \ 1 \right] 
\end{equation}
 is equal to $ 6\ln 2-2\ln(1+\sqrt{2})-\frac{16}{\pi}\operatorname{Im}\Li_2\left[(\sqrt{2}-1)i\right]$.
\end{theorem}

\begin{proof}
 Since
\begin{align*}
 \sum_{n\geq 1}\frac{\binom{4n}{2n}\binom{2n}{n}}{n\,64^n} 
 & = 2\int_{0}^{1}\left[\frac{2 \text{{\bf K}}\left(\sqrt{\frac{2y}{1+y}}\right)}{\sqrt{1+y}}-1\right]\frac{dy}{y} \\
 & = 4\int_{0}^{1}\left[\frac{\sqrt{2}}{\pi}\text{{\bf K}}(\sqrt{x})\sqrt{2-x}-1\right]\frac{dx}{x(2-x)} 
\end{align*}
\noindent we find that 
\begin{equation}\label{20180611305AM1A}
 \sum_{n\geq 1}\frac{\binom{4n}{2n}\binom{2n}{n}}{n\,64^n} 
 = \frac{\pi}{\sqrt{2}} + 
 6\ln 2-\frac{2\sqrt{2}}{\pi}\ln^2(1+\sqrt{2})-\frac{16}{\pi}G 
 - \frac{2}{\pi}\int_{0}^{1}\frac{\text{{\bf K}}(\sqrt{x})\,dx}{1+\sqrt{1-\frac{x}{2}}}. 
\end{equation}

 So, the whole problem of evaluating the series in \eqref{20180609433PM1A}
 boils down to finding the FL expansion of 
 $\frac{1} {1+\sqrt{1-\frac{x}{2}}}$:
\begin{align*}
 \frac{1} {1+\sqrt{1-\frac{x}{2}}} & = \sum_{n\geq 0}P_n(2x-1)\cdot 2\int_{1}^{+\infty}\frac{dv}{(v+1)^{3/2}\left(\sqrt{v}+\sqrt{v+1}\right)^{2n+1}}\\ 
 & = \sum_{n\geq 0}P_n(2x-1)\cdot 4\int_{0}^{(\sqrt{2}-1)^2}\frac{1-x}{(1+x)^2}x^n\,dx. 
\end{align*}
 Manipulate the above integral so as to obtain the following: 
\begin{align*}
 \sum_{n\geq 0}P_n(2x-1)\cdot 8 \left[ \frac{1}{2\sqrt{2}}(\sqrt{2}-1)^{2n+1} - \right. 
 \left. \int_{0}^{\sqrt{2}-1}(2n+1)x^{2n+1}\frac{dx}{1+x^2}\right].
\end{align*}

 From the FL expansion for $\text{{\bf K}}(\sqrt{x})$, we have that $$f(x)=\sum_{n\geq 0}P_n(2x-1) c_n \Longrightarrow 
 \int_{0}^{1}\text{{\bf K}}(\sqrt{x})\,f(x)\,dx = 2\sum_{n\geq 0}\frac{c_n}{(2n+1)^2},$$ letting $c_{n}$ denote a scalar 
 coefficient for $n \in \mathbb{N}_{0}$. By taking $f(x)=\frac{1}{1+\sqrt{1-\frac{x}{2}}}$ and considering the previous line, we 
 find that the definite integral in \eqref{20180611305AM1A} equals 
\begin{equation}\label{20180611307AM2A}
 4\sqrt{2}\sum_{n\geq 0}\frac{(\sqrt{2}-1)^{2n+1}}{(2n+1)^2}-16\int_{0}^{\sqrt{2}-1}\frac{\arctanh x}{1+x^2}\,dx. 
\end{equation}
\noindent We again encounter the series from \eqref{20190120951PM2A},
 and we thus find that 
 the expression in \eqref{20180611307AM2A} 
 must be equal to $$ \frac{\pi^2}{2\sqrt{2}}-\sqrt{2}\ln^2(1+\sqrt{2})-16\int_{0}^{\sqrt{2}-1}\frac{\arctanh x}{1+x^2}\,dx. $$ 
 The above integral may also be computed through the machinery of dilogarithms, yielding the following evaluation of 
 \eqref{20180611307AM2A}: $$ \frac{\pi^2}{2\sqrt{2}}-\sqrt{2}\ln^2(1+\sqrt{2})-8G+\pi\ln(1+\sqrt{2}) + 
 8\text{Im}\text{Li}_2\left[(\sqrt{2} - 1) i \right].$$ We thus obtain the desired result.
\end{proof}

\begin{theorem}
 The following equality holds: 
\begin{align*}
 & \sum _{n=0}^{\infty } \frac{\binom{2 n}{n}^2 (2 n+1)}{16^n (n+1)^4} \\ 
 & = 16-6 \pi ^2-32 \ln (2)+24 \ln ^2(2)+\frac{64 G-32 + 
 256 \cdot \text{\emph{Im}}\left(\text{\emph{Li}}_3\left(\frac{1 + 
 i}{2}\right)\right)}{\pi }.
\end{align*}
\end{theorem}

\newpage

\begin{proof}
 We have that 
 $$\sum_{n\geq 0}\frac{\binom{2n}{n}^2(2n+1)}{16^n(n+1)^4}
 = -\frac{4}{\pi}\int_{0}^{1}\frac{d}{dx}\text{{\bf E}}(\sqrt{x})\ln^2(x)\,dx,$$
 and we may evaluate 
 the right-hand side as 
\begin{align*}
 -\frac{4}{\pi}\Bigg[(2-\pi)+\sum_{n\geq 1}(-1)^n 
 & \Bigg(2\Bigg(\frac{2n+1}{n(n+1)}\Bigg)^2+4\frac{2n+1}{n(n+1)}H_{n-1}\Bigg) \\ 
 & \Bigg(-\frac{1}{2n+1}+\int_{0}^{1}\frac{2 x^{2n+2}}{1+x^2}\,dx\Bigg)\Bigg].
\end{align*}
 It is seen that the above expression 
 is equal to 
\begin{align*}
 = -\frac{4}{\pi}\big[ & -4\pi+8-16G+\frac{11}{12}\pi^3+(8\pi-16 G)\ln 2+\pi\ln^2 2 \\ 
 & - 32\,\text{Im }\text{Li}_3\big(\tfrac{1+i}{2}\big)-8\int_{0}^{1}\frac{\ln^2(1+x^2)}{1+x^2}\,dx\big].
\end{align*}
The desired result follows.
\end{proof}

\section{Generalized complete elliptic integrals}\label{20180614114PM2A}
 The study of the moments of $\text{{\bf K}}(\sqrt{x})$ and $\text{{\bf E}}(\sqrt{x})$ provides us with some new 
 results on generalized complete elliptic integrals. 
 While the following 
 definition is not new per se \cite{Takeuchi2016}, 
 the notation introduced below will be convenient.

\begin{definition}
\normalfont We define the generalized complete elliptic integral $ \mathfrak{J} $ over $[0,1)$ as 
\begin{align*} 
 \mathfrak{J}(x) 
 & = \int_{0}^{\pi/2}\left(\sqrt{1-x\sin^2\theta}\right)^3\,d\theta\\ 
 &= \frac{x-1}{3}\, \text{{\bf K}}(\sqrt{x}) + 
 \frac{4-2x}{3}\,\text{{\bf E}}(\sqrt{x}).
 \end{align*}
\end{definition}

 Inspired by the main results in the preceding sections, we are interested in:

\begin{itemize}
 \item the computation of the Taylor series for $ \mathfrak{J} $ at the origin; 
 \item the computation of the FL expansion for $ \mathfrak{J}$; and 
 \item the computation of its moments $\int_{0}^{1}\mathfrak{J}(x)\,x^\eta\,dx$.
\end{itemize}
 
 It is not difficult to see that $\mathfrak{J}(x) $ may be expanded as 
\begin{equation}\label{heartsuit1}
 \frac{3\pi}{2}\sum_{n\geq 0}\frac{\binom{2n}{n}^2}{16^n(1-2n)(3-2n)}x^n = \frac{\pi}{2}\cdot 
 {}_{2}F_{1}\!\!\left[ \begin{matrix} 
 -\tfrac{3}{2},\tfrac{1}{2}
 \vspace{1mm}\\ 
 1 \end{matrix} \ \Bigg| \ x \right], 
\end{equation}
\noindent and it is seen that the following equality holds: 
\begin{equation*} 
 \mathfrak{J}(x) = 48\sum_{n\geq 0}\frac{P_n(2x-1)}{(2n+5)(2n+3)(2n+1)(2n-1)(2n-3)}. 
\end{equation*}
 By multiplying the right-hand side of 
 \eqref{heartsuit1} by $x^\eta$ and performing a termwise integration we immediately have $$ \int_{0}^{1}\mathfrak{J}(x) 
 x^\eta\,dx = \frac{\pi}{2(1+\eta)}\cdot 
 {}_{3}F_{2}\!\!\left[ \begin{matrix} 
 -\tfrac{3}{2},\tfrac{1}{2},1+\eta
 \vspace{1mm}\\ 
 1,2+\eta \end{matrix} \ \Bigg| \ 1 \right]. $$
 From \eqref{201901211014PM2A}
 it follows that the $\eta^{\text{th}}$ moment of $ \mathfrak{J}(x)$ 
 is equal to 
 $$ 48\sum_{n\geq 0}(-1)^n \frac{\Gamma(n-\eta)\Gamma(1+\eta)}{\Gamma(-\eta)\Gamma(n+2+\eta)(2n + 5)(2n + 3)(2n 
 + 1)(2n - 1) (2n - 3)}. $$
 So, we find for $\eta > -1$ that 
\begin{align*}
 \int_{0}^{1}\mathfrak{J}(x) x^\eta\,dx 
 & = \frac{\pi}{2(\eta+1)}\cdot 
 {}_{3}F_{2}\!\!\left[ \begin{matrix} 
 -\tfrac{3}{2},\tfrac{1}{2},1+\eta
 \vspace{1mm}\\ 
 1,2+\eta \end{matrix} \ \Bigg| \ 1 \right] \\ 
 & = \frac{16}{15(\eta+1)}\cdot 
 {}_{3}F_{2}\!\!\left[ \begin{matrix} 
 -\tfrac{3}{2},1,-\eta
 \vspace{1mm}\\ 
 \tfrac{7}{2},2+\eta \end{matrix} \ \Bigg| \ -1 \right].
\end{align*}

\begin{corollary}
 The following equality holds for $\eta > -1$:
 $$ \frac{15\pi}{32} = \frac{ {}_{3}F_{2}\!\!\left[ \begin{matrix} -\tfrac{3}{2},1,-\eta \vspace{1mm}\\ 
 \tfrac{7}{2},2+\eta \end{matrix} \ \Bigg| \ -1 \right] }{ {}_{3}F_{2}\!\!\left[ \begin{matrix} -\tfrac{3}{2},\tfrac{1}{2},1+\eta 
 \vspace{1mm}\\ 1,2+\eta \end{matrix} \ \Bigg| \ 1 \right]}.$$ 
\end{corollary}

\begin{proof}
 This follows from 
 the expansions of $ \int_{0}^{1}\mathfrak{J}(x) x^\eta\,dx $ given above. 
\end{proof}

\begin{definition}
\normalfont The generalized complete elliptic integral $\mathfrak{J}_m$ 
 is given for all $x \in [0, 1)$ by 
 $$ \mathfrak{J}_m(x) = \int_{0}^{\pi/2}\left(1 - 
 x\sin^2\theta\right)^{m-\frac{1}{2}}\,d\theta. $$ 
\end{definition}

 We observe that 
 $\mathfrak{J}_{0}(x) = \text{{\bf K}}(\sqrt{x})$ and that $ \mathfrak{J}_1(x) = \text{{\bf E}}(\sqrt{x})$. Also, we have that 
 $\mathfrak{J}_2(x) = \mathfrak{J}(x)$. 

 By considering the FL expansions of $\mathfrak{J}_0,\mathfrak{J}_1,\mathfrak{J}_2$ 
 and by performing induction on $m$, we obtain many by-products. 

\begin{theorem}\label{tripleform}
 For any $\eta > -1$, the moment of $\mathfrak{J}_m(x)$
 of order $\eta$ 
 satisfies all of the properties 
\begin{align*} 
\int_{0}^{1}\mathfrak{J}_m(x)x^\eta\,dx & = 
 \frac{2\cdot 4^m}{\binom{2m}{m}(2m + 
 1)}\cdot 
 {}_{3}F_{2}\!\!\left[ \begin{matrix} 
-\eta,1,m+1
 \vspace{1mm}\\ 
\tfrac{3}{2},\tfrac{3}{2}+m \end{matrix} \ \Bigg| \ 1 \right] \\
 & = \frac{\pi}{2(1+\eta)}\cdot 
 {}_{3}F_{2}\!\!\left[ \begin{matrix} 
\tfrac{1}{2}-m,\tfrac{1}{2},1+\eta
 \vspace{1mm}\\ 
 1,2+\eta \end{matrix} \ \Bigg| \ 1 \right] \\ 
 & = \frac{2\cdot 4^m}{\binom{2m}{m}(2m+1)(\eta+1)}\cdot
 {}_{3}F_{2}\!\!\left[ \begin{matrix} 
 \tfrac{1}{2}-m,1,-\eta
 \vspace{1mm}\\ 
 \tfrac{3}{2}+m,2+\eta \end{matrix} \ \Bigg| \ -1 \right]. 
\end{align*}
\end{theorem}

\begin{proof}
 This follows by termwise integration
 of a function of the form $x^\eta \cdot \phantom{}_2 F_1$, and by 
 exploiting the ``trigonometric definition'' of $ \mathfrak{J}_m(x)$ along with Fubini's Theorem, 
 and by using FL expansions.
\end{proof}

\begin{corollary}\label{20180820419PM1A}
 The following equality holds for 
 $m \geq 0$ and $\eta > -1$: 
\begin{equation}\label{201808231137PM1A}
 \pi = \frac{4^{m+1}}{\binom{2m}{m}(2m+1)}\cdot 
 \frac{ {}_{3}F_{2}\!\!\left[ \begin{matrix} 
 \tfrac{1}{2}-m,1,-\eta
 \vspace{1mm}\\ 
 \tfrac{3}{2}+m,2+\eta \end{matrix} \ \Bigg| \ -1 \right]}{ 
 {}_{3}F_{2}\!\!\left[ \begin{matrix} 
\tfrac{1}{2}-m,\tfrac{1}{2},1+\eta
 \vspace{1mm}\\ 
 1,2+\eta \end{matrix} \ \Bigg| \ 1 \right]}.
\end{equation}
\end{corollary}

\begin{proof}
 This follows immediately from Theorem \ref{tripleform}.
\end{proof}

 By applying parameter derivatives to the equality \eqref{201808231137PM1A}, 
 we often end up 
 with results on harmonic summations for $\frac{1}{\pi}$. 
 For example, if we apply 
 the operator$\frac{\partial }{\partial n} \cdot \Big|_{n=0}$
 to both sides of \eqref{201808231137PM1A}
 and then let the parameter $\eta$ approach a natural number $j$, 
 we obtain an explicit evaluation of the infinite series 
$$ \sum _{i=0}^{\infty } \frac{ \binom{2 i}{i}^2 \left(2 H_{2 i}-H_i\right)}{16^{i} (i + j + 1)}
 = \frac{8 (j!)^2}{\pi} 
 \sum _{i=0}^j \frac{(2 i + 1) \left(H_{i-\frac{1}{2}}+\ln (2)\right) + 1}{(2 i + 1)^2 (j-i)! (i +j+1)!}.$$

 It seems worthwhile to consider the evaluation of integrals of the form 
 $$ \int_{0}^{1}x^{a}(1-x)^{b}\,\mathfrak{J}_m(x)^2\,dx \in \mathbb{Q} + \zeta(3)\,\mathbb{Q}$$
 for $m,a,b\in\mathbb{N}$. 
 For instance, in the case $m=0$ and $n=2$, from the FL expansion of $x(1-x)\text{{\bf K}}(\sqrt{x})$
$$ x(1-x) \text{{\bf K}}(\sqrt{x}) = 
 \sum_{n\geq 0}\frac{8(9-4n-4n^2)}{(2n+1)(4n^2+4n-15)^2}\,P_n(2x-1) $$
we have:
$$ \int_{0}^{1} x^2(1-x)^2\, \text{{\bf K}}^{2}(\sqrt{x}) \,dx = \frac{7}{2^{14}}\left(251\,\zeta(3)-18\right). $$
 The following formulas can also be obtained from our FL-based technique: 
\begin{align*}
 \int_{0}^{1}x\, \text{{\bf K}}^{2}(\sqrt{x}) \,dx 
 & =\tfrac{1}{4}\left(7\,\zeta(3)+2\right), \\ 
 \int_{0}^{1}x^2\, 
 \text{{\bf K}}^{2}(\sqrt{x}) \,dx 
 & = \tfrac{1}{64}\left(77\,\zeta(3)+34\right). 
\end{align*}
 
 Let us state in an explicit way the structure of the Taylor series and the FL expansions of our generalized complete elliptic 
 integrals, together with their special values at $x = \frac{1}{2}$. The following results may be proved by considering the 
 information about $\mathfrak{J}_0,\mathfrak{J}_1,\mathfrak{J}_2$ collected so far and by applying 
 induction on $m \in \mathbb{N}$: 
\begin{align*}
 \mathfrak{J}_m(x) & = \tfrac{\pi}{2}\cdot 
 {}_{2}F_{1}\!\!\left[ 
 \begin{matrix}
 -\tfrac{2m+1}{2},\tfrac{1}{2} \vspace{1mm}\\ 1
 \end{matrix} \ \Bigg| \ x \right] \\ 
 & = \frac{\pi}{2}(2m+1)!!\sum_{n\geq 0}\frac{\binom{2n}{n}^2\,x^n}{16^n (1-2n)\cdots (2m+1-2n)} \\ 
 & = 2(2m)!\sum_{n\geq 0}\frac{P_n(2x-1)}{(2m+1+2n)(2m-1+2n)\cdots(2m-1-2n)}, \\ 
 \mathfrak{J}_m\left(\tfrac{1}{2}\right) 
 & =2(2m)!\sum_{n\geq 0}\frac{\binom{2n}{n}(-1)^n}{4^n (2m+1+4n)(2m-1+4n)\cdots(2m-1-4n)} \\ 
 &= \frac{2\cdot 4^m}{\binom{2m}{m}(2m+1)}\cdot {}_{3}F_{2}\!\!\left[ 
 \begin{matrix}
 \tfrac{1-2m}{4},\tfrac{3-2m}{4},\tfrac{1}{2} \vspace{1mm}\\
\tfrac{2m+3}{4},\tfrac{2m+5}{4}
 \end{matrix} \ \Bigg| \ -1 \right] \\
 &= \int_{0}^{\pi/2}\left(1-\tfrac{1}{2}\sin^2\theta\right)^{\frac{2m-1}{2}}\,d\theta \\ 
 & = \frac{1}{2}\int_{0}^{1}\left(1-\frac{u}{2}\right)^{\frac{2m-1}{2}}\frac{du}{\sqrt{u(1-u)}}\\ 
 &= 2^{\frac{1-2m}{2}}\int_{0}^{1}\frac{(1+u^2)^{2m}\,du}{\sqrt{1-u^4}} \\ 
 & = e_m\, \text{{\bf E}}\left(\tfrac{1}{\sqrt{2}}\right) 
+ k_m\,\text{{\bf K}}\left(\tfrac{1}{\sqrt{2}}\right),\qquad e_m,k_m\in\mathbb{Q}. 
 \end{align*}

A generalization of Legendre's relation for the $\mathfrak{J}_m$ functions has already been proved by Shingo Takeuchi in \cite{Takeuchi2016}. 
By FL expansions, the computation of the integrals $\int_{0}^{1}\mathfrak{J}_m(x)^2\,dx$ boils down to the computation of the partial fraction decomposition of 
\begin{align*}
 & q_n^{(m)} = \frac{1}{(2n+2m+1)^2} \cdots 
 \frac{1}{(2n+3)^2} \cdot \frac{1}{(2n+1)^3} \cdots 
 \frac{1}{(2n-2m+1)^2}, 
\end{align*}
 which, regarded as a meromorphic function of the $n$ variable, 
 has a triple pole at $n=-\frac{1}{2}$ and double poles at $n\in\left\{-m-\frac{1}{2},\ldots,-\frac{3}{2},\frac{1}{2},\ldots, m-\frac{1}{2}\right\}$. We have 
$$ 2(2m)!\frac{(-1)^m (2n-2m-1)!!}{(2n+2m+1)!!}= \frac{2}{4^m}\sum_{k=0}^{2m}\binom{2m}{k}\frac{(-1)^k}{2n+2m+1-2k}$$
and by squaring both sides, multiplying them by $\frac{1}{2n+1}$ and applying $\sum_{n\geq 0}(\ldots)$, 
we obtain an evaluation of
 $ \int_{0}^{1}\mathfrak{J}_m(x)^2\,dx$ in terms of Ap\'{e}ry's constant $\zeta(3)$. 
 By the same argument, the integral $\int_{0}^{1}\mathfrak{J}_a(x)\mathfrak{J}_b(x)\,dx$ has a similar structure.

\section{Conclusion}\label{20180614120PM1A}
 We have investigated the relationships between classical hypergeometric sums and twisted 
 hypergeometric series, FL series, and the $\text{{\bf K}}$ and $\text{{\bf E}}$ mappings. 
 Our article may lead to new areas of research 
 into generalizing our main theorems. In particular, from the 
 results presented in Section \ref{20180614114PM2A}, we are especially interested in generalizing the $\mathfrak{J}_m$ 
 transformation as much as possible in a meaningful way, through a 
 suitable modification of our methods. 

\section*{References}

\bibliography{mybibfile}

\begin{thebibliography}{22}
\expandafter\ifx\csname natexlab\endcsname\relax\def\natexlab#1{#1}\fi
\expandafter\ifx\csname url\endcsname\relax
  \def\url#1{\texttt{#1}}\fi
\expandafter\ifx\csname urlprefix\endcsname\relax\def\urlprefix{URL }\fi

\bibitem[{Bauer(1859)}]{Bauer1859}
Bauer, G., 1859. Von den {C}oefficienten der {R}eihen von {K}ugelfunctionen
  einer {V}ariablen. J. Reine Angew. Math. 56, 101--121.

\bibitem[{Brafman(1951)}]{Brafman1951}
Brafman, F., 1951. Generating functions of {J}acobi and related polynomials.
  Proc. Amer. Math. Soc. 2, 942--949.

\bibitem[{Campbell(2018{\natexlab{a}})}]{Campbell2018New}
Campbell, J., 2018{\natexlab{a}}. New series involving harmonic numbers and
  squared central binomial coefficients.
\newline\urlprefix\url{https://hal.archives-ouvertes.fr/hal-01774708/}

\bibitem[{Campbell(2018{\natexlab{b}})}]{Campbell2018Ramanujan}
Campbell, J.~M., 2018{\natexlab{b}}. Ramanujan-like series for
  {$\frac{1}{\pi}$} involving harmonic numbers. Ramanujan J. 46~(2), 373--387.

\bibitem[{Campbell(2018{\natexlab{c}})}]{Campbell2018Series}
Campbell, J.~M., 2018{\natexlab{c}}. Series containing squared central binomial
  coefficients and alternating harmonic numbers.
\newline\urlprefix\url{https://hal.archives-ouvertes.fr/hal-01836014/}

\bibitem[{Campbell et~al.(2019)Campbell, D'Aurizio, and
  Sondow}]{CampbellDAurizioSondow}
Campbell, J.~M., D'Aurizio, J., Sondow, J., 2019. Hypergeometry of the
  parbelos. Amer. Math. Monthly.

\bibitem[{Campbell and Sofo(2017)}]{CampbellSofo2017}
Campbell, J.~M., Sofo, A., 2017. An integral transform related to series
  involving alternating harmonic numbers. Integral Transforms Spec. Funct.
  28~(7), 547--559.

\bibitem[{Chan et~al.(2013)Chan, Wan, and Zudilin}]{ChanWanZudilin2013}
Chan, H.~H., Wan, J., Zudilin, W., 2013. Legendre polynomials and
  {R}amanujan-type series for {$1/\pi$}. Israel J. Math. 194~(1), 183--207.

\bibitem[{Cohl and MacKenzie(2013)}]{MR3322254}
Cohl, H.~S., MacKenzie, C., 2013. Generalizations and specializations of
  generating functions for {J}acobi, {G}egenbauer, {C}hebyshev and {L}egendre
  polynomials with definite integrals. J. Class. Anal. 3~(1), 17--33.
\newline\urlprefix\url{https://doi.org/10.7153/jca-03-02}

\bibitem[{D'Aurizio and Trani(2018)}]{Daurizio2017}
D'Aurizio, J., Trani, S.~D., 2018. Surprising identities for the hypergeometric
  {$_4F_3$} function. Boll. Unione Mat. Ital. 11~(3), 403--409.
\newline\urlprefix\url{https://doi.org/10.1007/s40574-017-0142-0}

\bibitem[{El-Mikkawy and Cheon(2005)}]{ElCheon2005}
El-Mikkawy, M. E.~A., Cheon, G.-S., 2005. Combinatorial and hypergeometric
  identities via the {L}egendre polynomials---a computational approach. Appl.
  Math. Comput. 166~(1), 181--195.

\bibitem[{Erd\'{e}lyi et~al.(1953)Erd\'{e}lyi, Magnus, Oberhettinger, and
  Tricomi}]{MR0058756}
Erd\'{e}lyi, A., Magnus, W., Oberhettinger, F., Tricomi, F.~G., 1953. Higher
  transcendental functions. {V}ols. {I}, {II}. McGraw-Hill Book Company, Inc.,
  New York-Toronto-London, based, in part, on notes left by Harry Bateman.

\bibitem[{Gonz\'{a}lez(1954)}]{MR0064915}
Gonz\'{a}lez, M.~O., 1954. Elliptic integrals in terms of {L}egendre
  polynomials. Proc. Glasgow Math. Assoc. 2, 97--99.

\bibitem[{Holdeman(1970)}]{Holdeman1970}
Holdeman, Jr., J.~T., 1970. Legendre polynomial expansions of hypergeometric
  functions with applications. J. Mathematical Phys. 11, 114--117.

\bibitem[{Kirillov(1995)}]{MR1356515}
Kirillov, A.~N., 1995. Dilogarithm identities. Progr. Theoret. Phys.
  Suppl.~(118), 61--142, quantum field theory, integrable models and beyond
  (Kyoto, 1994).
\newline\urlprefix\url{https://doi.org/10.1143/PTPS.118.61}

\bibitem[{Levrie(2010)}]{Levrie2010}
Levrie, P., 2010. Using {F}ourier-{L}egendre expansions to derive series for
  {${1\over\pi}$} and {${1\over\pi^2}$}. Ramanujan J. 22~(2), 221--230.

\bibitem[{Takeuchi(2016)}]{Takeuchi2016}
Takeuchi, S., 2016. Legendre-type relations for generalized complete elliptic
  integrals. J. Class. Anal. 9~(1), 35--42.

\bibitem[{Wan and Zudilin(2012)}]{WanZudilin2012}
Wan, J., Zudilin, W., 2012. Generating functions of {L}egendre polynomials: a
  tribute to {F}red {B}rafman. J. Approx. Theory 164~(4), 488--503.

\bibitem[{Wan(2012)}]{Wan2012}
Wan, J.~G., 2012. Moments of products of elliptic integrals. Adv. in Appl.
  Math. 48~(1), 121--141.

\bibitem[{Wan(2013)}]{Wan2013}
Wan, J.~G., 2013. Random walks, elliptic integrals and related constant. Ph.D.
  thesis, University of Newcastle.

\bibitem[{Wan(2014)}]{Wan2014}
Wan, J.~G., 2014. Series for {$1/\pi$} using {L}egendre's relation. Integral
  Transforms Spec. Funct. 25~(1), 1--14.

\bibitem[{Zhou(2014)}]{Zhou2014}
Zhou, Y., 2014. Legendre functions, spherical rotations, and multiple elliptic
  integrals. Ramanujan J. 34~(3), 373--428.

\end{thebibliography}

\end{document}